\documentclass[11pt]{amsart}
\usepackage{amssymb,amsfonts, color,epsf}
\usepackage [cmtip,arrow]{xy}
\xyoption{all}
\usepackage {pb-diagram,pb-xy}
\usepackage{enumerate}
\usepackage{accents}
\usepackage[noautoscale]{youngtab}
\usepackage{subfigure}
\ifx\pdfpageheight\undefined
   \usepackage[dvips,colorlinks=true,linkcolor=blue,citecolor=red,%
      urlcolor=green]{hyperref}
   \usepackage[dvips]{graphicx}
   \makeatletter
   \edef\Gin@extensions{\Gin@extensions,.mps}
   \makeatother
\else
 \usepackage[pdftex]{graphicx}
   \usepackage[bookmarksopen=false,pdftex=true,breaklinks=true,%
      backref=page,pagebackref=true,plainpages=false,%
      hyperindex=true,pdfstartview=FitH,colorlinks=true,%
      pdfpagelabels=true,colorlinks=true,linkcolor=blue,%
      citecolor=red,urlcolor=green,hypertexnames=false%
      ]%
    {hyperref}
\fi
\usepackage{bm}

\usepackage{amsmath,amssymb,amsfonts}
\newtheorem{theorem}{Theorem}
\newtheorem*{theorem*}{Theorem}
\newtheorem{lemma}{Lemma}

\newtheorem{proposition}{Proposition}
\newtheorem*{proposition*}{Proposition}

\theoremstyle{definition}
\newtheorem{definition}{Definition}
\newtheorem{example}{Example}

\newtheorem{notation}{Notation}

\newtheorem{claim}{Claim}

\theoremstyle{remark}
\newtheorem{remark}{Remark}

\definecolor{DarkBlue}{rgb}{0,0.1,0.55}

\numberwithin{equation}{section}

\newcommand{\hide}[1]{}

\DeclareMathAlphabet{\mathpzc}{OT1}{pzc}{m}{it}
\newcommand{\Reeb}{\mathrm{Reeb}}

\usepackage{amssymb,color,epsf}
\usepackage{mathtools}
\usepackage{enumerate}

\usepackage{bm}
\usepackage{amsmath}
\usepackage{tabularx}
\usepackage{color, colortbl}
\usepackage{url}

\DeclareMathAlphabet{\mathpzc}{OT1}{pzc}{m}{it}
\usepackage[margin=1.5in]{geometry}

\usepackage [cmtip,arrow]{xy}
\xyoption{all}
\usepackage {pb-diagram,pb-xy}
\usepackage[mathscr]{eucal}
\usepackage[utf8]{inputenc}
\usepackage[T1]{fontenc}
\usepackage[english]{babel}
\usepackage{bm}
\usepackage{tabularx}

\usepackage{bm}

 \newcommand {\sign} {\mbox{\bf sign}}

\newcommand {\junk}[1]{}

\newcommand {\R} {\mathrm{R}}

\newcommand {\D}     {\mbox{\rm D}}

\newcommand {\Sphere}{\mbox{${\bf S}$}}     
\newcommand {\Disk}{\mbox{${\bf D}$}} 

\newcommand {\Q}         {\mathbb{Q}}

\newcommand {\ZZ} {\mathrm{Z}}
\newcommand {\RR} {{\mathcal R}}

\newcommand {\la}   {{\langle}}
\newcommand {\ra}   {{\rangle}}
\newcommand {\eps} {{\varepsilon}}
\newcommand {\E} {{\rm Ext}}

\newcommand {\PP}     {\mathbb{P}} 

\newcommand{\card}{\mathrm{card}}

\def\addots{\mathinner{\mkern1mu
\raise1pt\vbox{\kern7pt\hbox{.}}
\mkern2mu\raise4pt\hbox{.}\mkern2mu
\raise7pt\hbox{.}\mkern1mu}}

\newcommand{\HH}  {\mbox{\rm H}}

\DeclareMathOperator{\diag}{diag}

\newcommand{\x}{\mathbf{x}}

\newcommand{\y}{\mathbf{y}}

\newcommand{\supp}{\mathrm{Supp}}

\newcommand{\Cc}{\mathrm{Cc}}

\newcommand{\rR}{\mathrm{R}}

\newcommand{\nc}{\newcommand}
\newcommand{\rc}{\renewcommand}
\nc{\mc}{\mathcal}
\rc{\t}{\text}
\nc{\op}[1]{\operatorname{#1}}
\nc{\opcat}[1]{\mathbf{#1}}

\nc{\id}{\op{id}}

\nc{\umutnote}[1]{{\marginpar{\small \textcolor{blue}{#1}}}}

\nc{\cA}{\mc{A}}\nc{\cB}{\mc{B}}\nc{\cC}{\mc{C}}\nc{\cD}{\mc{D}}\nc{\cE}{\mc{E}}\nc{\cF}{\mc{F}}\nc{\cG}{\mc{G}}\nc{\cH}{\mc{H}}\nc{\cI}{\mc{I}}\nc{\cJ}{\mc{J}}\nc{\cK}{\mc{K}}\nc{\cL}{\mc{L}}\nc{\cM}{\mc{M}}\nc{\cN}{\mc{N}}\nc{\cO}{\mc{O}}\nc{\cP}{\mc{P}}\nc{\cQ}{\mc{Q}}\nc{\cR}{\mc{R}}\nc{\cS}{\mc{S}}\nc{\cT}{\mc{T}}\nc{\cU}{\mc{U}}\nc{\cV}{\mc{V}}\nc{\cW}{\mc{W}}\nc{\cX}{\mc{X}}\nc{\cY}{\mc{Y}}\nc{\cZ}{\mc{Z}}

\rc{\PP}{\mathbb{P}}
\rc{\AA}{\mathbb{A}}
\nc{\bbC}{\mathbb{C}}
\nc{\CC}{\mathbb{C}}

\nc{\code}[1]{{\texttt{#1}}}
\nc{\mcode}[1]{{\text{\texttt{#1}}}}
\nc{\xto}[1]{\raisebox{-0.03cm}{\scalebox{0.85}{$\,\xrightarrow{#1}\,$}}}
\nc{\xtonormal}[1]{\xrightarrow{#1}}
\nc{\xfrom}[1]{\xleftarrow{#1}}
\nc{\sidenote}[1]{\marginpar{\small #1}}

\nc{\Aff}{\opcat{Aff}}
\nc{\AffVar}{\opcat{AffVar}}
\nc{\ProjVar}{\opcat{ProjVar}}
\nc{\GAP}{\opcat{GrAlgPairs}}
\nc{\GA}{\opcat{GrAlg}}

\nc{\acc}{\mathrm{a.c.c}}
\nc{\GL}{\mathrm{GL}}

\nc{\Mod}{\t{-}\opcat{Mod}}
\nc{\Sub}{\opcat{Sub}}
\nc{\iso}{\cong}
\nc{\compose}{\circ}

\newcommand{\Rc}{\mathrm{R}}

\newcommand{\SIGN}{\mathrm{SIGN}}

\newcommand{\crit}{\mathrm{crit}}
\newcommand{\Crit}{\mathrm{Crit}}

\newcommand\restr[2]{{
  \left.\kern-\nulldelimiterspace 
  #1 
  \vphantom{\big|} 
  \right|_{#2} 
  }}

\title[On the Reeb spaces of definable maps]
{
On the Reeb spaces of definable maps
}
\author{Saugata Basu}
\address{Department of Mathematics,
Purdue University, West Lafayette, IN 47906, U.S.A.}
\email{sbasu@math.purdue.edu}

\author{Nathanael Cox}
\address{Department of Mathematics,
Purdue University, West Lafayette, IN 47906, U.S.A.}
\email{cox175@math.purdue.edu}

\author{Sarah Percival}
\address{Department of Mathematics,
Purdue University, West Lafayette, IN 47906, U.S.A.}
\email{sperciva@math.purdue.edu}

\date{\textbf{\today}}
\thanks{
Basu was partially supported by NSF grants
CCF-1618918 and DMS-1620271. Cox and Percival were partially supported by NSF grant DMS-1620271.
 }

\begin{document}
\begin{abstract}
We prove that the Reeb space of a proper definable map $f:X \rightarrow Y$  in an arbitrary o-minimal expansion of a real closed field is realizable as a proper definable quotient. This result can be seen as
an o-minimal analog of Stein factorization of proper morphisms in algebraic geometry.
We also show that the Betti numbers of the Reeb space of $f$ can be arbitrarily large compared to those of $X$, unlike in the special case of Reeb graphs of manifolds.
Nevertheless, in the special case when $f:X \rightarrow Y$ is a semi-algebraic
map and $X$ is closed and bounded, we prove a singly exponential upper bound on the Betti numbers of  the Reeb space of $f$ in terms of the number and degrees of the polynomials defining $X,Y$, and $f$. 
\end{abstract}

\maketitle
\tableofcontents

\section{Introduction}
Given a topological space $X$ and a continuous function $f \colon X \to \mathbb{R}$, define an equivalence relation $\sim$ on $X$ by setting $x \sim x'$ if and only if $f(x) = f(x')$ and $x$ and $x'$ are in the same connected component of $f^{-1}(f(x)) =  f^{-1}(f(x'))$. The space $X/\sim$ is called the \emph{Reeb graph} of $f$, denoted $\Reeb(f)$. The concept of the Reeb graph was introduced by Georges Reeb in \cite{Reeb} as a tool in Morse theory. The notion of the Reeb graph can be  generalized to the notion of  \emph{Reeb space} by letting $f \colon X \to Y$, where $Y$ is any topological space. 
Burlet and de Rham first introduced the Reeb space in \cite{Rham} as the \emph{Stein factorization} of a map $f$, but their work was limited to bivariate, generic, smooth mappings. 
Existence of Stein factorization for more general morphisms in algebraic geometry is proved in \cite[III, Corollary 11.5]{Hartshorne},  
and is closely related to  the well-known Zariski's Main Theorem \cite[III, Corollary 11.4]{Hartshorne} (see Remarks~\ref{rem:Stein} and~\ref{rem:Stein2}  for the connection between Stein factorization in algebraic geometry and the results of the current paper).

From the point of view of applied topology, Reeb spaces have been investigated from both a theoretical and practical perspective.
Edelsbrunner et al.\ defined the Reeb space of a multivariate piecewise linear mapping on a combinatorial manifold in \cite{RSP}, and they proved results regarding the local and global structure of such spaces. Expanding on this work, Patel \cite{Patel} produced an algorithm to construct the Reeb space of a mapping $f$. Mapper, introduced in \cite{Singh}, gives a discrete approximation of the Reeb space of a multivariate mapping; this allows for more efficient computation of the underlying data structure. Munch et al.\ \cite{Munch} define the \emph{interleaving distance} for Reeb spaces to show the convergence between the Reeb space and Mapper.

In this paper, we investigate Reeb spaces from the point of view of topological complexity. Our motivation is to understand how
topologically complicated the Reeb space of a map can become in terms of the complexity of the map itself. In order to obtain meaningful results we restrict ourselves to the category of maps \emph{definable in an o-minimal expansion of  
a real closed field
$\R$} (for example, one can take $\R = \mathbb{R}$)
and
in particular, to \emph{semi-algebraic} maps
(see Section~\ref{sec:om} for a quick overview of o-minimality).  

\begin{remark}
\label{rem:Rspaces}
We remark here that in \cite{de-Silva-et-al} Reeb graphs are considered for so-called constructible $\mathbb{R}$-spaces
that the authors introduce. Compact definable sets in any o-minimal expansion of $\mathbb{R}$ are 
examples of constructible $\mathbb{R}$-spaces. Our choice of studying Reeb spaces within the framework of o-minimal geometry stems from the fact that o-minimality is now a widely accepted framework for studying
tame geometry and also allows us to prove effective upper bounds on the topology of Reeb spaces in special situations of interest -- for instance, for Reeb spaces of semi-algebraic maps.
\end{remark}

The notion of o-minimal structures has its origins in model theory but has since become a widely accepted framework for studying
``tame geometry''. 
The definable sets and maps of an o-minimal structure satisfy many uniform finiteness properties (similar to those of semi-algebraic sets) 
while allowing much richer families of sets and maps.
We refer the reader to the survey by Wilkie \cite{Wilkie2009} for the origin and motivation of this notion of tameness. The reader will also find many applications of interest. 

Our first result is that the Reeb spaces of ``tame'' maps are themselves tame. More precisely, we prove that the quotient map corresponding to the Reeb space of a proper definable map can be realized as 
a proper definable map (Theorem~\ref{thm:definability-of-Reeb} below). This implies as a special case that the Reeb spaces of proper semi-algebraic maps can be realized as semi-algebraic quotients. Theorem~\ref{thm:definability-of-Reeb} can be viewed as the definable analog of the theorem~\cite[III, Corollary 11.5]{Hartshorne} on the existence of Stein factorization for proper morphisms in
algebraic geometry (see Remark~\ref{rem:Stein} below).
Another significance of this result is that it makes it possible to ask for an algorithm to
semi-algebraically describe this semi-algebraic quotient using results from the well developed area of algorithmic semi-algebraic geometry \cite{BPRbook2}. 
A naive approach would be to  make the proof of Proposition~\ref{prop:definability-of-quotient} in \cite{Dries} (in the semi-algebraic case) algorithmic. However,
the proof of Proposition~\ref{prop:definability-of-quotient} makes heavy use of definable (in this case semi-algebraic) triangulations. 
Algorithms with the best known complexity for computing semi-algebraic  triangulations  (see for instance \cite[Chapter 5]{BPRbook2}) use the technique of \emph{cylindrical algebraic decomposition} which has  intrinsically a doubly exponential complexity.  So this approach is unlikely to yield anything better than an algorithm with a doubly exponential complexity. 
However, it must be noted that the problem of computing a description of the Reeb space is much more special
than the problem of computing the description of quotients by an arbitrary semi-algebraic equivalence relations,
since the equivalence relation in the case of the Reeb space is very geometric. 
Many problems in algorithmic semi-algebraic geometry with a geometric or topological flavor-- such as the problem of 
computing semi-algebraic descriptions of the semi-algebraically connected components of a given semi-algebraic set, computing certain topological invariants like the Euler-Poincar\'e characteristics etc. admit algorithms  with singly exponential complexity (see for examples \cite[Chapters 13-16]{BPRbook2}). 
They rely on a common technique called the \emph{critical point method}.

In this paper we prove a singly exponential upper bound on the Betti numbers of the Reeb space
of a semi-algebraic map.  Driven by the analogy with the problems mentioned above, this singly exponential upper bound suggests that an algorithm having a singly exponential 
complexity should exist for computing a description of the Reeb space using the critical point method.
However, the technique introduced in this paper is not sufficient for this purpose, since the upper bound that
we prove comes from approximating the cohomology of the Reeb space by the $E_1$-term of a spectral
sequence converging to it. But this is not sufficient to describe the Reeb space itself. Doing so would require
additional techniques and 
we do not pursue this question further in this paper, leaving it for future work.

It is known \cite[page 141]{Edelsbrunner-Harer} that the sum of the Betti numbers of the \emph{Reeb graph} of a map 
$f:X \rightarrow \R$ is bounded from above by the sum of the Betti numbers of $X$. 
We show that this is false for more general maps
by exhibiting a couple of natural examples of sequences of maps  $(f_n : X_n \rightarrow Y_n)_{n > 0}$, 
such that the sum of the Betti numbers of the Reeb space of $f_n$ is arbitrarily
large compared to that of $X_n$. In view of these examples, it makes sense to ask whether it is still possible to bound the
Betti numbers of the Reeb space of a map $f$ in terms of some measure of the  ``complexity'' of the map $f$.
In particular, if the map is semi-algebraic, then one can measure the complexity of the map  by the number and degrees of the
polynomials defining the map. We are then led to
the problem of studying the topological complexity of Reeb spaces of semi-algebraic maps.

While studying the topological complexity of Reeb spaces of semi-algebraic maps is a natural mathematical question
on its own, 
another motivation is related to the algorithmic question mentioned earlier concerning the design of efficient
algorithms for computing a semi-algebraic description of the Reeb space of a 
semi-algebraic map. It is a meta theorem in algorithmic semi-algebraic geometry that upper bounds on topological complexity of objects are closely related to the worst-case complexity of algorithms computing the topological invariants of such
objects.  Thus, a singly exponential upper bound on the Betti numbers of the Reeb space of a semi-algebraic map opens up the possibility of being able to compute the Betti numbers of the Reeb space. The singly exponential upper bound on the Betti numbers of the Reeb space of a semi-algebraic map may also hint that one could compute a semi-algebraic description of the Reeb space with an algorithm having a singly exponential complexity bound.

The problem of bounding the topological complexity (for example measured in terms of Betti numbers or the number of homotopy types of fibers) of semi-algebraic sets or maps in terms of the parameters of the formula defining them has a long history (see \cite{BPR10} for a survey). Bounds on these quantities which are doubly exponential in the dimension or the number of variables usually follow from the fact that semi-algebraic sets admit semi-algebraic triangulations of at most  doubly
exponential size. Singly exponential upper bounds are more difficult and usually involve more careful arguments involving
Morse inequalities and other inequalities coming from certain spectral sequences \cite{OP,T,Milnor2,BPR02,GaV,BPRbook2}. To the best of our knowledge, the problem of bounding the Betti numbers of the Reeb space of a semi-algebraic map has not been 
 considered before. In this
 paper we prove a singly exponential upper bound on the Betti numbers of the Reeb space of a semi-algebraic map
 $f:X \rightarrow Y$, where $X$ is a closed and bounded semi-algebraic set, in terms of the number and the degrees of the
 polynomials defining $X,Y$ and $f$ (cf.\ Theorem~\ref{thm:main} below). 
 
The rest of the paper is organized as follows: 
In Section~\ref{sec:om}, we recall the basic definitions related to o-minimality. In Section~\ref{sec:om-Reeb}, we prove the definability of Reeb spaces of proper definable maps. In Section~\ref{sec:example}, we describe examples showing
that the Betti numbers of the Reeb space of a definable map $f:X \rightarrow Y$ can be arbitrarily large compared to those of $X$.
We also give a proof of the inequality $b_1(\Reeb(f))\leq b_1(X)$ for definable proper maps $f:X \rightarrow Y$ with $X$ connected,
using a spectral sequence that plays an important role in this paper (this inequality was proved previously using alternative techniques by
Dey et al.\ \cite{Dey-et-al-2017}).
Finally, in Section~\ref{sec:quantitative}, we prove a singly exponential upper bound on the sum of the Betti numbers of the Reeb space of a 
proper semi-algebraic map in terms of the number and degrees of the polynomials defining the map. We defer
the proof of a key intermediate result (Proposition~\ref{prop:BElim}) to the Appendix (Section~\ref{sec:appendix})
since it is long and technical.
We conclude in Section \ref{sec:conclusion} with some open problems.

\section{Basic definitions}
\label{sec:om}
We first recall the important  model theoretic notion of o-minimality which plays an important role in what follows.
  
\subsubsection{O-minimal Structures}
O-minimal structures were invented and first studied by
Pillay and Steinhorn in the pioneering papers
\cite{PS1,PS2}, motivated by the prior work of van den Dries \cite{Dries5}. 
Later, the theory was further
developed through contributions of other researchers, most notably
van den Dries, Wilkie, Rolin, and Speissegger, amongst others
\cite{Dries2,Dries3,Dries4,Wilkie,Wilkie2,Rolin}. We particularly
recommend the book by van den Dries \cite{Dries} and the notes by
Coste \cite{Michel2} for an easy introduction to the topic as well as for the
proofs of the basic results that we use in this paper.

\begin{definition}[o-minimal structure]
\label{def:o-minimal}
An o-minimal structure over 
a real closed field 
$\Rc$  (or equivalently an o-minimal expansion of $\Rc$) is a sequence 
${\mathcal S}(\Rc) = ({\mathcal S}_n)_{n \in {\mathbb N}}$ where 
each ${\mathcal S}_n$ is a collection of subsets of $\Rc^n$
(called the {\em definable sets} in the structure) satisfying the 
following axioms (following the exposition in \cite{Michel2}): 

\begin{enumerate}
\item
All algebraic subsets of $\Rc^n$ are in ${\mathcal S}_n$.
\item
The class ${\mathcal S}_n$ is closed under complementation and
finite unions and intersections.
\item
If $A \in {\mathcal S}_m$ and $B \in {\mathcal S}_n$ then
$A \times B \in {\mathcal S}_{m+n}$.
\item
If $\pi: \Rc^{n+1} \rightarrow \Rc^{n}$ is the projection map on the
first $n$ coordinates and $A \in {\mathcal S}_{n+1}$, then 
$\pi(A)  \in {\mathcal S}_n$.
\item
The elements of ${\mathcal S}_1$ are finite unions of points
and intervals.
(Note that these are precisely the subsets of $\Rc$ which are
definable by a first-order formula in the language of the reals 
with one free variable.) 
\end{enumerate}
A map $f:X \rightarrow Y$ between two definable sets $X$ and $Y$ 
is \emph{definable}
if its graph is 
a definable set.
Note that for any definable map $f:X \rightarrow Y$, there exists a finite partition $(X_i)_{i\in I}$ of $X$ into definable subsets such that
$f$ restricted to each $X_i$ is continuous. In light of this, for rest of this paper we use the term ``definable map'' to mean a map that is definable and continuous.
\end{definition}

The class of semi-algebraic subsets of $\R^n, n >0$, where $\R$ is a real closed field,
is one obvious example of an o-minimal structure,
but in fact there are much richer classes of sets which have been proven
to be o-minimal. The class of \emph{sub-analytic sets} is one such example \cite{Wilkie2}.


We now consider quotients by definable equivalence relations.

\begin{definition}
Let $E \subset X \times X$ be a definable equivalence relation on a definable set $X$. A \emph{definable quotient} of $X$ by $E$ is a pair
$(p,Y)$ consisting of a definable set $Y$ and a definable surjective map $p:X \rightarrow Y$ such that
\begin{enumerate}[(i)]
\item
$(x_1,x_2) \in E  \Leftrightarrow p(x_1) = p(x_2)$ for all $x_1,x_2 \in X$;
\item
$p$ is definably identifying; that is, for all definable $K \subset Y$, if $p^{-1}(K)$ is closed in $X$, then $K$ is closed in $Y$.
\end{enumerate}
We say that the definable quotient $(p,Y)$ is \emph{definably proper} if $p$ is a definably proper map.
This means that if $X \subset \R^m$, and $Y \subset \R^n$, then for every definable subset $K\subset Y$ with $K$ closed and bounded in $\Rc^n$, 
$p^{-1}(K)\subset X$ is closed and bounded in $\Rc^m$.
\end{definition}

\begin{definition}
A definable equivalence relation $E \subset X \times X$ is said to be \emph{definably proper} if the two maps
$\op{pr}_1, \op{pr}_2:E \rightarrow X$ are definably proper.
\end{definition}

We will use the following proposition which appears in \cite{Dries}:

\begin{proposition}\cite[page 166]{Dries}
\label{prop:definability-of-quotient}
Let $X$ be a definable set and  $E \subset X \times X$ a definably proper equivalence relation on $X$. Then $X/E$ exists as a definably proper quotient of $X$.
\end{proposition}

\section{The Reeb space of a definable map $f:X \rightarrow Y$}
\label{sec:om-Reeb}
We now fix an o-minimal expansion of  a real closed field $\R$.
Let $X \subset \R^n$ be a closed and bounded definable set, and $f:X \rightarrow Y$ be a definable map.

\begin{definition}
The Reeb space of the map $f$, henceforth denoted $\Reeb(f)$, is the topological space
$X/{\sim}$, equipped with the quotient topology, where $x \sim x'$ if and only if $f(x) = f(x')$, and $x,x'$ belong to the same connected component of $f^{-1}(f(x))$. 
(Here the topology on $X$ is the definable topology whose
open sets are the definable open subsets of $X$).
\end{definition}

\begin{remark}
\label{rem:connectivity}
Note that a definable (resp.\ semi-algebraic) set $S \subset \R^k$ is connected if and only if $S$ is definably (resp.\ semi-algebraically) path-connected, i.e.
for all $x,y \in S$, there exists a definable (resp.\ semi-algebraic) path $\gamma:[0,1] \rightarrow S$ with $\gamma(0) = x, \gamma(1) = y$.
\end{remark}

\begin{theorem}
\label{thm:definability-of-Reeb}
Let $X \subset \R^n$ be a closed and bounded definable set, and $f:X \rightarrow Y$ be a definable map.
Then, the space $\Reeb(f) \triangleq X/{\sim}$ 
is a definably proper quotient. 
In other words, let $X \subset R^n$ be a closed and bounded definable set, and $f\colon X \to Y$ be a definable map. Then
there exists a definable set $Z$, and a  proper definable map $\psi:X \rightarrow Z$ 
and a homeomorphism $\theta: \Reeb(f) \rightarrow Z$ such that the following diagram commutes:
\[
\xymatrix{
&X \ar[ld]^\phi \ar[rd]_\psi& \\
\Reeb(f)=X/{\sim} \ar[rr]^\theta && Z
}
\] 
(here $\phi$ is the quotient map).
In particular, $\Reeb(f)$ is homeomorphic to a definable set.
\end{theorem}

\begin{remark}
\label{rem:properness}
The assumption of 
that $X$ is closed and bounded 
is needed. For example, suppose $X =\R^2 \setminus \mathbf{0}$ and $f:X \rightarrow \R$ is the projection map forgetting the second coordinate. Then, each fiber $f^{-1}(x)$ has one connected component if $x \neq 0$,
and $f^{-1}(0)$ has two connected components. The Reeb space of $f$ is homeomorphic to  
$\R$
with a doubled point, and 
is not definable.
\end{remark}

\begin{remark}
\label{rem:Stein}
Theorem~\ref{thm:definability-of-Reeb} can also be seen as a definable analog of Stein factorization for projective morphisms ~\cite[III, Corollary 11.5]{Hartshorne} which states that
``every projective morphism $f:X \rightarrow Y$ of Noetherian schemes factors as $f = g \circ f'$,  with $g:Y' \rightarrow Y$ a finite morphism, and $f': X \rightarrow Y'$ a projective morphism with connected fibers.''

Note that this result is valid for Noetherian schemes over fields of any characteristic and the connectivity of the fibers is  with respect to the
the Zariski topology of the underlying topological space of the corresponding scheme.

Notice also that the scheme $Y'$ plays the role of Reeb space of $f$. More precisely, the underlying 
topological space of the scheme $Y'$ is the Reeb space of the map induced by  the morphism $f$ between the underlying topological spaces of the schemes $X$ and $Y$.

We discuss below a simple illustrative example assuming  only a minimum familiarity with some concepts from algebraic geometry, such as affine and projective  schemes over a field $K$, homogeneous coordinates etc. (the reader can consult \cite{Hartshorne} for definitions).
\begin{example}
\label{eg:Stein}
Let $K$ be any field. Recall that a morphism $f:X \rightarrow Y$ between $K$-schemes is  projective
if it factors as $f = g \circ f'$, where $f': X \rightarrow \mathbb{P}^N_K \times Y$ is a closed immersion
for some $N > 0$, and $g: \mathbb{P}^N_K \times Y \rightarrow Y$ is the projection 
\cite[page 103]{Hartshorne}.

Consider in $\mathbb{A}^3_K$ the subscheme consisting of the union two lines defined by the equations parametrized by $T$:
$X_3=T, X_2=0$ and  $X_3=0, X_1=0$.
Then for $T\neq 0$, the two lines are disjoint, but at $T=0$ the subscheme is connected.
In order to produce an example of a projective morphism we need to projectivize the picture.

We will use $(X_0:X_1:X_2:X_3)$ to denote the homogeneous coordinates of $\mathbb{P}_K^3$.
Let  $X$ denote the subscheme of $\mathbb{P}^3_K \times \mathbb{A}^1_K$ defined by 
intersection of the ideals $(X_2,X_3-X_0T)$ and  $(X_1,X_3)$.
The fibers of $X$ with respect to the projection to $\mathbb{A}^1_K$ consist generically of two skew
lines in $\mathbb{P}^3_K$ except over the point defined by $T=0$. Over $T=0$, the two lines meet, and so the fiber is connected.

Let 
$Y' \subset \mathbb{P}^1_K \times \mathbb{A}^1_K$ be the subscheme
defined by the polynomial
\[
W_1(W_1 - W_0 T)
\]
and let $g: Y' \rightarrow \mathbb{P}^1_K = Y$ be the restriction to $Y'$ of the projection morphism
$\mathbb{P}^1_K \times \mathbb{A}^1_K \rightarrow \mathbb{A}^1_K$. Note that the fibers of $g$ consist of two points except over $T=0$ where it is a
doubled point.

Let 
\[
f':   X \rightarrow    Y'
\]
be the morphism defined by
\[
((X_0:X_1:X_2:X_3), T) \mapsto ((X_0:X_3),T).
\]

Note that $X_0$ and $X_3$ cannot simultaneously be $0$ on $X$ and so $f'$ is well defined.
Note  also that $f'$ factors through the closed embedding
\[
f'': X \rightarrow \mathbb{P}^3_K \times Y'
\]
defined by 
\[
((X_0:X_1:X_2:X_3), T) \mapsto ((X_0:X_1:X_2:X_3),(X_0:X_3),T),
\]
followed by the projection to $Y'$. Hence, $f'$ is a projective morphism. 

It can now be verified that:
\begin{enumerate}[1.]
\item
$f = g \circ f' $, 
\item$g$ is a finite morphism, and 
\item
$f'$ has connected fibers.
\end{enumerate}

\end{example}
\end{remark}

We will now prove Theorem~\ref{thm:definability-of-Reeb}.
\begin{proof}[Proof of Theorem~\ref{thm:definability-of-Reeb}]
We first claim that the relation, 
``$x \sim x'$ if and only if $f(x) = f(x')$, and $x,x'$ belong to the same connected component of $f^{-1}(f(x))$''
is a definably proper equivalence relation. Using Hardt's triviality theorem for o-minimal structures \cite{Dries,Michel2},  we have that there exists a finite 
definable partition of $Y$ into locally closed definable sets $(Y_\alpha)_{\alpha \in I}$, $y_\alpha \in Y_\alpha$, and definable homeomorphisms $\phi_\alpha: Y_\alpha \times f^{-1}(y_\alpha) \rightarrow f^{-1}(Y_\alpha)$ 
such that the following diagram commutes for each $\alpha \in I$:

\[
\xymatrix{
Y_\alpha \times f^{-1}(y_\alpha) \ar[rd]_{\pi_1} \ar[rr]^{\phi_\alpha}  && f^{-1}(Y_\alpha)\ar[ld]^{f|_{f^{-1}(Y_\alpha)}} \\
&Y_\alpha&
}
\]
(here $\pi_1$  is the projection to the first factor in the direct product). 
For each ${\alpha \in I}$, let $(C_{\alpha,\beta})_{\beta \in J_\alpha}$ be the connected components of $f^{-1}(y_\alpha)$, 
and for each $\alpha \in I, \beta\in J_{\alpha}$, let $D_{\alpha,\beta} = \phi_\alpha(Y_\alpha \times C_{\alpha,\beta})$.

Let
\[
E = \bigcup_{\alpha \in I,\beta \in J_\alpha} (\phi_\alpha \times \phi_\alpha)((Y_\alpha \times C_{\alpha,\beta}) \times_{\pi_1} (Y_\alpha \times C_{\alpha,\beta})),
\]
where $(Y_\alpha \times C_{\alpha,\beta}) \times_{\pi_1} (Y_\alpha \times C_{\alpha,\beta})$ is the definable subset of 
$(Y_\alpha \times f^{-1}(y_\alpha)) \times (Y_\alpha \times f^{-1}(y_\alpha)) $
defined by 
\[
((y,x),(y',x')) \in (Y_\alpha \times C_{\alpha,\beta}) \times_{\pi_1} (Y_\alpha \times C_{\alpha,\beta}) \Leftrightarrow y = y', x,x' \in C_{\alpha,\beta}.
\]
It is clear that $E$ is a definable subset of $X \times X$, and that $x \sim x'$ if and only if $(x,x') \in E$.

Since $X$ is assumed to be closed and bounded,  if we can show that $E$ is closed in $X \times X$,
it would follow that 
$E$ is a definably proper equivalence relation, and we can apply Proposition~\ref{prop:definability-of-quotient}.

The rest of the proof is devoted to showing that $E$ is a closed definable subset of $X \times X$.
For each $\alpha \in I, \beta \in J_\alpha$, let 
\[
E_{\alpha,\beta} = (\phi_\alpha \times \phi_\alpha)((Y_\alpha \times C_{\alpha,\beta}) \times_{\pi_1} (Y_\alpha \times C_{\alpha,\beta})).
\]

Since $E = \bigcup_{\alpha \in I,\beta \in J_\alpha} E_{\alpha,\beta}$, in order to prove that $E$ is closed it suffices to prove that for each $\alpha \in I, \beta \in J_\alpha$,  
\[
\overline{E_{\alpha,\beta}} \subset E,
\]
where $\overline{E_{\alpha,\beta}}$ is the closure of $E_{\alpha,\beta}$ in $X \times X$.

It follows from the curve selection lemma for o-minimal structures \cite{Michel2} that for every $z \in \overline{E_{\alpha,\beta}}$
there exists a definable curve $\gamma:[0,1] \rightarrow E_{\alpha,\beta}$ with $\gamma(0) = z$, $\gamma((0,1]) \subset E_{\alpha,\beta}$. Thus, in order to prove that $\overline{E_{\alpha,\beta}} \subset E$, it suffices to show that for each definable curve 
$\gamma: (0,1] \rightarrow E_{\alpha,\beta}$, $z_0 = \lim_{t \rightarrow 0} \gamma(t)  \in E$.

Let $\gamma: (0,1] \rightarrow E_{\alpha,\beta}$ be a definable curve, and suppose that $\lim_{t\rightarrow 0} \gamma(0) \not\in E_{\alpha,\beta}$. 
Otherwise, $\lim_{t\rightarrow 0} \gamma(0) \in E_{\alpha,\beta} \subset E$, and we are done.

For $t \in (0,1]$, let $y_t = f(\gamma(t))$
and let $(x_t,x_t') \in (\phi_\alpha \times \phi_\alpha)((Y_{\alpha} \times C_{\alpha,\beta}) \times_{\pi_1} (Y_{\alpha} \times C_{\alpha,\beta}))$ be such that $\gamma(t) = (x_t,x_t')$. Note that 
$f(x_t) = f(x_t') = y_t$. 
Finally, let $z_0 = (x_0,x_0') = \lim_{t \rightarrow 0} \gamma(t)$.

Since, $z_0 \not\in E_{\alpha,\beta}$ by assumption and $\gamma((0,1]) \subset E_{\alpha,\beta}$, 
there exists $t_0 > 0$ such that
$\lambda = f \circ \gamma|_{(0,t_0]}: (0,t_0] \rightarrow Y_\alpha$ is an injective definable map
and $\lim_{t \rightarrow 0}\lambda(t) = y_0
= f(x_0) = f(x_0') \in Y_{\alpha'}$ for some $\alpha' \in I$.
We need to show that $x_0$ and $x_0'$ belong to the same connected component of $f^{-1}(y_0)$, which would imply that
$(x_0,x_0') \in E$.

Let $D_{\alpha,\beta,\gamma} = f^{-1}(\lambda((0,t_0])) \cap D_{\alpha,\beta}$ and let $g:D_{\alpha,\beta,\gamma} \rightarrow (0,t_0]$ be defined by
$g(x) = \lambda^{-1}(f(x))$ (which is well defined by the injectivity of $\lambda$). Note that for each $t \in (0,t_0]$, $g^{-1}(t)$ is 
definably homeomorphic to 
$C_{\alpha,\beta}$, 
and hence is connected. It also follows from Hardt's triviality theorem that there exists $t_0' \in (0,t_0]$ and a definable homeomorphism
$\theta:  g^{-1}(t_0') \times (0,t_0'] \rightarrow g^{-1}((0,t_0'])$ such that the following diagram commutes:
\[
\xymatrix{
g^{-1}(t_0') \times (0,t_0'] \ar[rr]^\theta \ar[rd]^{\pi_2} &&g^{-1}((0,t_0']) \ar[ld]^g \\
&(0,t_0']&
}
\]

Extend $\theta$ continuously to a definable map 
$\overline{\theta}:  g^{-1}(t_0') \times [0,t_0] \rightarrow \overline{g^{-1}((0,t_0'])}$ by setting $\overline{\theta}(x,0) = \lim_{t \rightarrow 0} \theta(x,t)$.
Finally, let 
$\theta': g^{-1}(t_0') \rightarrow f^{-1}(y_0)$ be the definable map obtained by setting 
$\theta'(x) = \overline{\theta}(x,0)$.

Note that since $g^{-1}(t_0')$ is connected, 
$\theta'(g^{-1}(t_0'))$ is connected as well, since it is the image of a connected set under a continuous map. 
Also note that for each $t \in (0,t_0']$, we have that $x_t,x'_t \in D_{\alpha,\beta,\gamma}$ and $f(x,t) = f(x'_t) = \lambda(t)$, 
hence $x_t,x'_t \in g^{-1}(t)$, and thus
$x_0,x_0' \in \theta'(g^{-1}(t_0'))$.
Moreover, $f(x_0) = f(x_0') = y_0$. Therefore, since $\theta'(g^{-1}(t_0'))$ is connected, $x_0$ and $x_0'$ belong to the same connected component of $f^{-1}(y_0)$.

This shows that $(x_0,x_0') \in E$, which in turn implies that $E$ is closed in $X \times X$. 

The fact that $\Reeb(f)$ exists as a definably proper quotient now follows from Proposition~\ref{prop:definability-of-quotient}.
\end{proof}

\begin{remark}
\label{rem:algorithm}
Theorem~\ref{thm:definability-of-Reeb} opens up an algorithmic problem of actually realizing the Reeb space as a definable
quotient in the special case where the o-minimal structure is that of semi-algebraic sets and maps. More precisely, the problem
is to design an algorithm that, given a proper semi-algebraic map $f:X \rightarrow Y$, will compute a description of  a semi-algebraic map $g:X \rightarrow Z$ such that following diagram commutes
\[
\xymatrix{
&X \ar[ld]^\phi \ar[rd]_g& \\
\Reeb(f)=X/{\sim} \ar[rr]^\theta && Z
}
\] 
(here $\phi$ is the quotient map and $\theta$ is a homeomorphism).
The complexity of the algorithm will then depend on the number and degrees of the polynomials defining $X,Y$,
and the graph of the semi-algebraic map $f$. In this paper, we do not pursue this algorithmic problem any further leaving it for  future work. 
\end{remark}

\section{The Betti numbers of the Reeb space of $f \colon X \to Y$ can exceed that of $X$}
\label{sec:example}

\begin{definition}
\label{def:Betti}
Any closed and bounded definable set is finitely triangulable \cite[Theorem 4.4]{Michel2}.
For any closed and bounded definable set $X$ and $i \geq 0$,  let $b_i(X)$ denote the $i$-th Betti number (that is, the dimension of the 
$i$-th simplicial homology group of a definable triangulation of $X$ with coefficients in $\mathbb{Q}$), and we denote $b(X) = \sum_{i \geq 0} b_i(X)$.
It is well known that the Betti numbers so defined do not depend on the chosen triangulation of $X$.
\end{definition}

In  \cite[page 141]{Edelsbrunner-Harer} it is noted that the  inequality
$b(\Reeb(f)) \leq b(X)$ holds 
for arbitrary maps $f:X \rightarrow \R$.

We first show that the same is not true for Reeb spaces of more general maps by giving several examples.

\begin{example}
\label{eg:counter-example}
Consider the closed $n$-dimensional disk $\Disk^n$ with $n \geq 1$, and let ${\sim}$ be the equivalence relation identifying all points on the boundary of $\Disk^n$. Then $\Disk^n/\sim \; \cong \Sphere^n$, where $\Sphere^n$ is the $n$-dimensional sphere. Let $f_n$ denote the quotient map $f_n: \Disk^n \to \Sphere^n$. The fibers of $f_n$ consist of either one point or the boundary $\Sphere^{n-1}$ of $\Disk^n$, and hence $\Reeb(f_n) \cong \Sphere^n$ for all $n > 1$. Note that $b_0(\Disk^n) = 1$ and $b_i(\Disk^n) = 0$ for all $i > 0$. Moreover, $b_0(\Sphere^n) = 1$, $b_n(\Sphere^n) = 1$, and $b_i(\Sphere^n)=0, i\neq 0,n$.
Thus, we have for $n > 1$,
\begin{eqnarray*}
b(\Disk^n) &=&1, \\
b(\Reeb(f_n)) &=& 2.
\end{eqnarray*}
More generally, for $k \geq 0$, let
\[ f_{n,k} = \underbrace{f \times \cdots \times f}_{\mbox{$k$ times}}: 
\underbrace{\Disk^n \times \cdots \times \Disk^n}_{\mbox{$k$ times}} \longrightarrow 
\underbrace{\Sphere^{n}\times \cdots \times \Sphere^{n}}_{\mbox{$k$ times}}.
\]
Using the same argument as before,
for $n >1$ and $k > 0$,  
\[\Reeb(f_{n,k}) \cong  \underbrace{\Sphere^{n}\times \cdots \times \Sphere^{n}}_{\mbox{$k$ times}}.\]
Thus, 
\begin{eqnarray*}
b_0(\underbrace{\Disk^n \times \cdots \times \Disk^n}_{\mbox{$k$ times}}) &=& 1,\\
b_i((\underbrace{\Disk^n \times \cdots \times \Disk^n}_{\mbox{$k$ times}})) &=& 0, \: i > 0, \\
\end{eqnarray*}
and hence
\[b(\underbrace{\Disk^n \times \cdots \times \Disk^n}_{\mbox{$k$ times}}) = 1.\]
Moreover, for $n > 1$, 
\begin{eqnarray*}
b_i(\Reeb(f_{n,k})) &=& 0 \mbox{ if $n\not| \;i$ or if $i > nk$},  \\
b_i(\Reeb(f_{n,k})) &=& \binom{k}{i/n} \mbox{ otherwise}, \\
\end{eqnarray*}
and hence for $n > 1$,
\[b(\Reeb(f_{n,k})) = 2^k.\]

This example shows  that even for definably proper maps $f:X \rightarrow Y$,
the individual as well as the total Betti numbers of $\Reeb(f)$ can  be arbitrarily large compared to those of $X$.
\end{example}

Our second example comes from the topology of compact Lie groups, in particular the complex unitary group:

\begin{example}
\label{sec:counter-example3}
For $n >0$,  let $U(n)$ denote the group of $n \times n$ complex unitary matrices, and let $T^n \subset U(n)$ denote the maximal torus. (Note that $T^n$ is the group of $n \times n$ unitary diagonal matrices
$\diag(z_1,\ldots,z_n)$ with $|z_i| = 1, 1 \leq i \leq n$, and is thus homeomorphic to the product of $n$ circles.)  
Denote the quotient map by $\pi_n: U(n) \rightarrow U(n)/T^n$.
We have that:
\begin{eqnarray}
\label{eqn:unitary1}
b(U(n)/T^n) &=& n! \mbox{ (see \cite[Theorem 4.6]{Mimura-Toda-book})}, \\
\label{eqn:unitary2}
b(U(n)) &=&  2^n \mbox{ (see \cite[Corollary 3.11]{Mimura-Toda-book})}.
\end{eqnarray}
Observing that the fibers of $\pi_n$ are all connected,
one has that $\Reeb(\pi_n) \cong U(n)/T^n$, and it follows from \eqref{eqn:unitary1} and \eqref{eqn:unitary2} that for all $n \geq 4$,
\[
b(\Reeb(\pi_n)) = n!  \geq 2^n = b(U(n)).
\]
\end{example}

We note that recently Dey et al.\ \cite{Dey-et-al-2017} have shown that

\begin{equation}
\label{eqn:b1}
b_1(\Reeb(f)) \leq b_1(X)
\end{equation}
if $f:X \rightarrow Y$ is a proper map and $X$ is connected.  
Notice that the examples given above do not violate this bound 
since the stated inequalities involve only the sum rather than the individual Betti numbers.

We sketch below an alternative proof of the inequality \eqref{eqn:b1}
for a proper definable map $f:X \rightarrow Y$, with $X$ connected,  using an inequality coming from a 
spectral sequence associated to the quotient map 
$\phi: X \rightarrow \Reeb(f)$.
This spectral sequence also plays a key role in the proof of the main result in this paper.

More precisely, for a proper definable surjective map $g: A \rightarrow B$, Gabrielov, Vorobjov and Zell  \cite{GVZ04} proved that there 
exists a spectral sequence (which we write as a cohomological spectral sequence for convenience) which converges to $\HH^*(B)$.
This spectral sequence is referred to as the \emph{descent spectral sequence} of $g$ below and 
its $E_1$-term  is given by
\[
E_1^{p,q} = \HH^q(\underbrace{A \times_g\cdots \times_g A}_{p+1}),
\] 
where 
\[
\underbrace{A \times_g\cdots \times_g A}_{p+1} = \{(a_0,\ldots,a_p) \in A^{p+1} \mid g(a_0) = \cdots = g(a_p) \}.
\]
\begin{theorem}\cite{Dey-et-al-2017}
\label{thm:Dey}
Let  $X$ be a connected definable set, and $f:X \rightarrow Y$ a proper map definable map.

Then, 
\begin{equation*}
b_1(\Reeb(f)) \leq b_1(X)
\end{equation*}
\end{theorem}

\begin{proof}[Alternative proof of Theorem~\ref{thm:Dey} using a spectral sequence argument]
First note that if $X$ is connected, then so is $\underbrace{X \times_\phi \cdots \times_\phi X}_{p+1}$,
and $\dim(E_1^{p,0}) = 1$ for all $p \geq 0$. Moreover, the differential $d^{p,0}_1:E_1^{p,0} \rightarrow E_1^{p+1,0}$ has rank 
$0$ or $1$ depending on whether $p$ is even or odd, respectively.
This implies that $E_2^{p,0} = 0$ for all $p > 0$ in the descent spectral sequence of  the quotient map $\phi:X \rightarrow \Reeb(f)$.
Moreover, notice that $E_1^{0,1} \cong \HH^1(X)$, and hence
\[
\dim(E_1^{0,1}) = b_1(X).
\]

Since the spectral sequence converges to $\HH^{p+q}(\Reeb(f))$,  the following inequality holds for each $n \geq 0$ and $r \geq 1$:

\begin{eqnarray}
\label{eqn:spectral1}
\HH^{n}(\Reeb(f)) &\leq& \sum_{p+q=n} \dim(E_r^{p,q}).
\end{eqnarray}

Moreover, for $r\geq r'$ and for any $p,q$,
\begin{eqnarray}
\label{eqn:spectral2}
\dim(E_r^{p,q}) &\leq & \dim(E_{r'}^{p,q}),
\end{eqnarray}
since $E_r^{p,q}$ is a sub-quotient of $E_{r'}^{p,q}$.

It follows from the inequalities \eqref{eqn:spectral1} and \eqref{eqn:spectral2} with $n=1$, $r'=1$, and $r=2$, that
\begin{eqnarray*}
b_1(\Reeb(f)) &\leq & \dim(E_2^{0,1}) + \dim(E_2^{1,0}) \\
&\leq& \dim(E_1^{0,1}) + \dim(E_2^{1,0}) \\
&=& b_1(X) + 0 \\
&=& b_1(X).
\end{eqnarray*}
\end{proof}

\begin{remark}
We note here that  an inequality  (cf.\ inequality \eqref{eqn:descent}) coming from the consideration of the $E_1$-term of the spectral sequence  of the map $\phi$ 
plays a 
key role in the proof of Theorem~\ref{thm:main}, which is the main result of this paper.
\end{remark}

\section{Quantitative Bounds}
\label{sec:quantitative}
We now consider the problem of bounding effectively from above the Betti numbers of the Reeb space of a definable 
continuous map.  
We have seen from Example~\ref{eg:counter-example}  that, given a continuous semi-algebraic map
$f:X \rightarrow Y$, $b(\Reeb(f))$ can be arbitrarily large compared to $b(X)$, unlike in the case of Reeb graphs (i.e. when
$\dim(Y) \leq 1$). In this section, we prove an upper bound on $b(\Reeb(f))$ in terms of the ``semi-algebraic'' complexity of the
map $f$. 

We first introduce some more notation.

\begin{notation}
  \label{not:sign-condition} 
  Let $\rR$ be a real closed field.
  For any finite family of polynomials $\mathcal{P} \subset \rR [ X_{1} , \ldots ,X_{k} ]$, w
  e call an element 
  $\sigma \in \{0,1,-1 \}^{\mathcal{P}}$ a \emph{sign condition} on $\mathcal{P}$. 
  For
  any semi-algebraic set $Z \subset \rR^{k}$ and sign condition $\sigma \in
  \{ 0,1,-1 \}^{\mathcal{P}}$, we denote by $\RR (\sigma ,Z)$ the
  semi-algebraic set defined by
  \[
  \{ \mathbf{x} \in Z \mid \sign (P (
  \mathbf{x})) = \sigma (P),  P \in \mathcal{P} \},
  \]
  and call it the
  \emph{realization} of $\sigma$ on $Z$. 
  
  More generally, we call any 
  Boolean formula $\Phi$ with atoms $P \{ =,>,< \} 0, P \in \mathcal{P}$, a \emph{$\mathcal{P}$-formula}. 
  We call the realization of $\Phi$,
  namely the semi-algebraic set
  \begin{eqnarray*}
    \RR (\Phi , \rR^{k}) & = & \{ \mathbf{x} \in \rR^{k} \mid
    \Phi (\mathbf{x}) \},
  \end{eqnarray*}
  a \emph{$\mathcal{P}$-semi-algebraic set}. Finally, we call a Boolean
  formula without negations and with atoms $P \{\geq, \leq \} 0$ (resp. $P \{>, < \} 0$), $P\in \mathcal{P}$, a 
  \emph{$\mathcal{P}$-closed} (resp. \emph{$\mathcal{P}$-open}) formula, and we call
  the realization, $\RR(\Phi , \rR^{k})$, a \emph{$\mathcal{P}$-closed} (resp. a  a \emph{$\mathcal{P}$-open})
  semi-algebraic set.
  
  We will denote by $\SIGN(\mathcal{P})$ the set of \emph{realizable sign conditions of $\mathcal{P}$}, i.e.
  \[
  \SIGN(\mathcal{P}) = \{ \sigma \in \{0,1,-1\}^{\mathcal{P}} \mid \RR(\sigma,\rR^k) \neq \emptyset \}.
  \]
  
Finally, for any semi-algebraic  set $S$, we will denote the set of its semi-algebraically 
connected components by $\Cc(S)$.
\end{notation}

We prove the following theorem.

\begin{theorem}
\label{thm:main}
Let $S \subset \R^n$ be a bounded $\mathcal{P}$-closed semi-algebraic set, and $f=(f_1,\ldots,f_m): S \rightarrow \R^m$ 
be a polynomial map. Suppose that $s = \card(\mathcal{P})$ and the maximum of the degrees of the polynomials in $\mathcal{P}$
and $f_1,\ldots,f_m$ is bounded by $d$. Then,
\[
b(\Reeb(f)) \leq (s d)^{(n+m)^{O(1)}}.
\] 
\end{theorem}

The rest of the paper is devoted to the proof of Theorem~\ref{thm:main}. We first outline the main idea behind the proof.

\subsection{Outline of the proof of Theorem~\ref{thm:main}}
\label{subsec:outline}
We first replace the map $f:S \rightarrow \R^m$, by a new map $\tilde{f}:\tilde{S} \rightarrow \R^m$,  where
$\tilde{S} \subset \R^n \times \R^m$ and $\tilde{f}$ is the restriction to $\tilde{S}$ of the projection map to $\R^m$, such that the following diagram commutes:
\[
\xymatrix{
\tilde{S} \: \ar[rd]_{\tilde{f}} \ar@{^{(}->}[rr]  && \R^n \times \R^m \ar[ld]^{{\pi_m}} \\
&\R^m&
}
\]
From the
definitions it is evident that $\Reeb(f)$ and $\Reeb(\tilde{f})$ are homeomorphic. We next prove that there exists a semi-algebraic partition of $\R^m$ of controlled complexity (more precisely given by the connected components of the realizable sign conditions of a family of polynomials of singly exponentially bounded degrees and cardinality) into connected semi-algebraic sets 
$C$, such that the connected components of the fibers $\tilde{f}^{-1}(z)$ are in 1-1 correspondence with each other
as $z$ varies over $C$. Moreover each of these connected components $C$ is described by
a quantifier-free first order formula and the complexity of these formulas  (i.e. the number of polynomials appearing in the formula and their respective degrees) is bounded singly exponentially 
(see Proposition~\ref{prop:parametrized-cc} below for the precise formulation of this statement).

The proof  of  this result (Proposition~\ref{prop:parametrized-cc})  
uses an intermediate result -- namely 
Proposition~\ref{prop:BElim} -- whose long and somewhat technical proof is deferred to the Appendix (Section~\ref{sec:appendix}).
The fact that the semi-algebraically connected components of a semi-algebraic set can be described efficiently (with singly exponential complexity) is a consequence of  a result in \cite{BPRbook2} (Proposition~\ref{prop:cc} below).

Next, we use the fact that 
the canonical surjection $\phi:\tilde{S} \rightarrow \Reeb(\tilde{f})$ is a proper semi-algebraic map. We then use an
inequality proved in \cite{GVZ04} (see Proposition~\ref{prop:descent} below) to obtain an upper bound on the Betti numbers of the image of a proper semi-algebraic map $F:X\rightarrow Y$
in terms of the sum of the Betti numbers of various fiber products 
$X \times_F \cdots \times_F X$  of the same map. Recall that for $p \geq 0$, the $(p+1)$-fold fiber product is given by
\[
\underbrace{ X  \times_F \cdots \times_F X}_{\mbox{$(p+1)$-times}}  \triangleq \{(x^{(0)},\ldots,x^{(p)}) \in X^{p+1} \;\mid\; F(x^{(0)}) =\cdots =F(x^{(p)})\}.
\]
 
Proposition~\ref{prop:parametrized-cc} provides us with a well controlled description  (i.e. by quantifier-free first order formulas
involving singly exponentially any polynomials of singly exponentially bounded degrees)  of the fibered products
$\tilde{S} \times_{\tilde{f}} \cdots \times_{\tilde{f}}\tilde{S}$.
Finally, using these descriptions and results on bounding the Betti numbers of general semi-algebraic sets in terms of the number and degrees of polynomials defining them (cf.\  Proposition~\ref{prop:GV07} below) we obtain the claimed bound
on $\Reeb(f)$.

In order to make the above summary precise we first need to state some preliminary results.

\subsection{Parametrized description of  connected components}
\label{subsec:prelim}

The following proposition, which states that given any finite family of polynomials 
\[
\mathcal{P} \subset \R[X_1,\ldots,X_k,Y_1,\ldots,Y_\ell],
\]
where $\R$ is a real closed field,
there exists a semi-algebraic partition of $\R^\ell$ of controlled complexity which
has good properties with respect to $\mathcal{P}$, will play a crucial role in the proof of Theorem~\ref{thm:main}.

We will use the following notation in the rest of the paper.
\begin{notation}
\label{not:pi}
We will denote by
$\pi_{Y}:\rR^{k+\ell}   \rightarrow \rR^\ell$ 
the projection to the last $\ell$ (denoted by  $Y= (Y_1,\ldots,Y_\ell)$) coordinates.
For any semi-algebraic 
subset $S \subset \rR^{k+\ell}$  and $T \subset \R^\ell$, we will denote by
$S_T= S \cap \pi_{Y}^{-1}(T)$. If $T = \{\y\}$, we will write $S_\y$ in stead of $S_{\{\y\}}$.
\end{notation}
 
\begin{proposition}
\label{prop:parametrized-cc}
Let $\R$ be a real closed field, and let $\mathcal{P} \subset \R[X_1,\ldots,X_k,Y_1,\ldots,Y_\ell]$ be a finite set of polynomials of degrees bounded by $d$, with $\card(\mathcal{P}) =s$. Let $S \subset \R^k \times \R^\ell$ be a $\mathcal{P}$-semi-algebraic set. Then
there exists a finite set of polynomials $\mathcal{Q} \subset \R[Y_1,\ldots,Y_\ell]$
such that $\card(\mathcal{Q})$ and the degrees of polynomials in $\mathcal{Q}$
are bounded by $(sd)^{(k+\ell)^{O(1)}}$, and $\mathcal{Q}$ has the following additional property.

For each 
$\sigma \in  \SIGN(\mathcal{Q}) \subset  \{0,1,-1\}^{\mathcal{Q}}$ 
and 
$C \in \Cc(\RR(\sigma,\R^\ell))$,
there exists
\begin{enumerate}[(i)]
\item
an index set $I_{\sigma,C}$,
\item
 a finite family of polynomials $\mathcal{P}_{\sigma,C}  \subset \R[X_1,\ldots,X_k,Y_1,\ldots,Y_\ell]$,
and 
\item
$\mathcal{P}_{\sigma,C}$-formulas, $(\Theta_\alpha(\overline{X},\overline{Y}))_{\alpha \in I_{\sigma,C}}$, 
\end{enumerate}
such that
\begin{enumerate}
\item $\Theta_\alpha(\x,\y) \Rightarrow \y \in C$;\\
\item for each $\y \in C$ 
and 
each $D \in \Cc(S_C)$,
there exists a unique $\alpha \in I_{\sigma,C}$ such that 
$\RR(\Theta_\alpha(\cdot,\y)) = D_\y$ and 
$D_\y \in \Cc(S_\y)$. 
\end{enumerate}
\end{proposition}

The proof of Proposition~\ref{prop:parametrized-cc} will use the following result on efficient descriptions of the connected components
of semi-algebraic sets which can easily be deduced from \cite[Theorem 16.3]{BPRbook2} and which we state without proof.

\begin{proposition}
\label{prop:cc}
  \label{16:the:dcc}
  Let $\R$ be a real closed field and
  let $\mathcal{P} = \{P_{1} , \ldots ,P_{s} \} \subset \R
  [X_{1} , \ldots ,X_{k} ]$ with $\deg (P_{i} ) \leq d$ for $1 \leq i
  \leq s$ and let a semi-algebraic set $S$ be defined by a $\mathcal{P}$-formula. 
  Then there exists an algorithm that outputs
  quantifier-free semi-algebraic descriptions of all the 
  connected components of $S$. The number of polynomials  that appear in the output is bounded by
  $s^{k+1} d^{O (k^{4} )}$, while the degrees of the polynomials  are bounded by $d^{O (k^{3} )}$. 
\end{proposition}

In order to prove Proposition~\ref{prop:parametrized-cc} we will also need the following proposition.

\begin{proposition}
\label{prop:BElim}
Let $\rR$ be a real closed field.
Let $\mathcal{P} \subset \rR[X_1,\ldots,X_k,Y_1,\ldots,Y_\ell]$ be a finite set of polynomials with degrees bounded by $d$ and with $\card(\mathcal{P}) =s$, and let $S \subset \rR^{k} \times \rR^\ell$ be a bounded $\mathcal{P}$-closed
semi-algebraic set. Then
there exists a finite set of polynomials, $\mathcal{Q} \subset \R[Y_1,\ldots,Y_\ell]$,
with 
\begin{eqnarray*}
\card(\mathcal{Q}) &\leq& (sd)^{(k+\ell)^{O(1)}}, \\
\max_{Q \in \mathcal{Q}} \deg(Q) &\leq& d^{(k+\ell)^{O(1)}},
\end{eqnarray*} 
such that  for  each  
$\sigma \in  \SIGN(\mathcal{Q})$, 
$C \in \Cc(\RR(\sigma,\R^\ell))$,
$\y \in C$, and $D \in \Cc(S_C)$,
$D_\y \in \Cc(S_\y)$.
\end{proposition}

\begin{proof}
See Section~\ref{sec:appendix} (Appendix).
\end{proof}

\begin{remark}
\label{rem:sheaf}
In the case $\R = \mathbb{R}$, Proposition~\ref{prop:BElim} can be deduced 
from a somewhat more general theorem \cite[Theorem 4.21]{Basu-sheaf} 
about constructible sheaves. However, to  explain this properly one would need to recall
basic definitions and functorial properties of constructible sheaves which would put an
unreasonable burden on the reader.  So we prefer to give 
a more self-contained proof avoiding the sheaf-theoretic formalism which can be found in the Appendix
(Section~\ref{sec:appendix}).

The connection of Proposition~\ref{prop:BElim}  
with the theory of constructible sheaves is the following. We assume that the reader is familiar with familiar with basic definitions sheaf theory 
(which can be found for example in \cite{Godement}).
Let $\mathcal{F}$ be the sheaf on $\R^\ell$ defined as the sheafification of the presheaf that associates to each open subset $U \subset \R^\ell$, the finite-dimensional $\Q$-vector space $\HH^0(\pi_Y^{-1}(U \cap S,\Q))$. For each $\y \in \R^\ell$,
the stalk $\mathcal{F}_\y$ of the sheaf $\mathcal{F}$ at $\y$ is then isomorphic to 
$\HH^0(S_\y,\Q)$. The property of the finite set of polynomials $\mathcal{Q}$ guaranteed by  Proposition~\ref{prop:BElim}  -- namely, that for  each 
$\sigma \in  \SIGN(\mathcal{Q}) \subset  \{0,1,-1\}^{\mathcal{Q}}$, 
$C \in \Cc(\RR(\sigma,\R^\ell))$,
$\y \in C$, and 
for $D \in \Cc(C_S)$,  
$D_\y \in \Cc(S_\y)$ --
can be expressed in the language of sheaves as follows.

The set $\mathcal{Q}$ induces a semi-algebraic partition of $\R^\ell$ into a finite number of locally closed
semi-algebraic sets, $\RR(\sigma,\R^\ell), \sigma \in \SIGN(\mathcal{Q})$, such that for each
$\sigma \in \SIGN(\mathcal{Q})$, the stalks $\mathcal{F}_\y$ are locally constant. This implies in particular that the
sheaf $\mathcal{F}$ is a constructible sheaf. 

Let $j$ denote the inclusion $j:S \hookrightarrow \R^{k+\ell}$, and $\Q_S = j_* j^{*}\Q_{\R^{k+\ell}}$, where 
$\Q_{\R^{k+\ell}}$ denotes the constant sheaf with stalks isomorphic to $\Q$ on $\R^{k+1\ell}$.
Then, 
\[
\mathcal{F} \cong R^0\pi_{Y,*} \Q_S
\] 
(i.e. $\mathcal{F}$ is isomorphic to the zero-th higher direct image of the constant sheaf supported on $S$ under the map $\pi_Y$). 
The general theorem (\cite[Theorem 4.21]{Basu-sheaf})
alluded to in the first paragraph of this remark gives a singly exponential bound on the ``complexity'' of the higher direct images of a constructible sheaf in terms of the complexity of the sheaf itself. The complexity of a constructible sheaf is defined in terms of the complexity of the smallest semi-algebraic partition on which the stalks are locally constant. Proposition~\ref{prop:BElim} could in principle be deduced from  \cite[Theorem 4.21]{Basu-sheaf},
since the partition $\R^{k+\ell}$ into the semi-algebraic set $S$ and its complement has the property that the
stalks of $\Q_S$ are locally constant on the sets of this partition, and applying the theorem only for the 
$0$-th higher direct image $\Q_S$  under the (proper) map $\pi_Y$.
\end{remark}

We are now in a position to prove Proposition~\ref{prop:parametrized-cc}.

\begin{proof}[Proof of Proposition~\ref{prop:parametrized-cc}]
Let $\Phi(\overline{X},\overline{Y})$ be the $\mathcal{P}$-closed formula describing $S$.

First apply Proposition~\ref{prop:BElim} to obtain a set of polynomials $\mathcal{Q} \subset \R[Y_1,\ldots,Y_\ell]$ 
with degrees and cardinality bounded by $(sd)^{(k+\ell)^{O(1)}}$, and for  each connected component $C$ of each realizable sign condition 
$\sigma \in  \SIGN(\mathcal{Q}) \subset  \{0,1,-1\}^{\mathcal{Q}}$, 
each $y \in C$, and 
for each connected component of $D$ of $\pi_Y^{-1}(C) \cap S$, 
$D_y = \pi_Y^{-1}(y) \cap D$ is a  connected component of $S_y = \pi_Y^{-1}(y) \cap S$.

Next using  Proposition~\ref{prop:cc} obtain for each realizable sign condition $\sigma$ of $\mathcal{Q}$, and for each connected component of $C$ of $\RR(\sigma,\R^\ell$), a quantifier-free formula $\Phi_{\sigma,C}(\overline{Y})$ describing $C$.

Now using Proposition~\ref{prop:cc} one more time, obtain for each $\sigma,C$, and each connected component 
$D_\alpha$ of the semi-algebraic set defined by $\Phi_{\sigma,C}(\overline{Y}) \wedge \Phi(\overline{X},\overline{Y})$, a quantifier-free formula 
$\Theta_{\alpha}(\overline{X},\overline{Y})$ describing $D_\alpha$. 
\end{proof}

\subsection{Bounding the topology of the image of a polynomial map}
\label{subsec:descent}

The following proposition proved in \cite{GVZ04} allows one to bound the Betti numbers of the image of a
closed and bounded definable set $X$ under a definable map $F$ 
in terms of the Betti numbers of  the iterated fibered product of $X$ over $F$. More precisely:
 
\begin{proposition} \cite{GVZ04}
\label{prop:descent}
Let $F: X \rightarrow Y$ be a definable continuous map, and $X$ a closed and bounded definable set. Then, for for all $p \geq 0$,
\[
b_p(F(X)) \leq \sum_{\substack{i,j \geq 0\\i+j = p}} b_i(\underbrace{X \times_F \cdots \times_F X}_{(j+1)}).
\]
\end{proposition}

\subsection{Bounds on the Betti numbers of semi-algebraic sets}
Finally, in order to prove Theorem~\ref{thm:main}, we will need singly exponential upper bounds on the Betti numbers 
of semi-algebraic sets in terms of the number and degrees of the polynomials appearing in any quantifier-free formula defining the set. We give a brief overview of these results. The key result that we will need in the proof of Theorem~\ref{thm:main} is  Proposition~\ref{prop:GV07}.
(We refer the reader to \cite[Chapter 6]{BPRbook2} for the definition of homology groups of semi-algebraic
sets which are only locally closed and not necessarily closed and bounded.)

\subsubsection{General Bounds}
\label{subsubsec:general}
The first results on bounding the Betti numbers of real varieties were proved by
Ole{\u\i}nik and Petrovski{\u\i} \cite{OP}, Thom \cite{T}, and Milnor \cite{Milnor2}.
Using a Morse-theoretic argument and Bezout's theorem they proved the following proposition which appears in \cite{BPR02} and makes more precise an earlier result which appeared in
\cite{B99} :
\begin{proposition} \cite{BPR02}
\label{prop:B99}
If $S \subset \R^k$ is a $\mathcal{P}$-closed semi-algebraic set, then
\begin{eqnarray}
\label{eqn:B99}
b(S) \leq \sum_{i=0}^{k} \sum_{j=0}^{k-i} \binom{s+1}{j} 6^j d(2d-1)^{k-1},
\end{eqnarray}
where $s = \card(\mathcal{P}) > 0$ and $d = \max_{P \in \mathcal{P}} \deg(P)$.
\end{proposition}

Using an additional ingredient (namely, a technique to replace an arbitrary semi-algebraic set by a locally
closed one with a very controlled increase in the number of polynomials used to describe the given set),
Gabrielov and Vorobjov \cite{GaV} extended Proposition~\ref{prop:B99} to arbitrary 
$\mathcal{P}$-semi-algebraic sets with only a small increase in the bound. Their result 
in conjunction with Proposition~\ref{prop:B99} gives the following proposition.

\begin{proposition} \cite{GV07,BPRbook2}
\label{prop:GV07}
If $S \subset \R^k$ is a $\mathcal{P}$-semi-algebraic set, then
\begin{eqnarray}
\label{eqn:GaV}
b(S) \leq \sum_{i=0}^{k} \sum_{j=0}^{k-i} \binom{2ks+1}{j} 6^j d(2d-1)^{k-1},
\end{eqnarray}
where $s = \card(\mathcal{P})$ and $d = \max_{P \in \mathcal{P}} \deg(P)$.
\end{proposition}

We will also use the following bound on the number of connected components of the realizations of 
all realizable sign conditions of a family of polynomials proved in \cite{BPR02}.
\begin{proposition}
  \label{7:the:theorem1}
Let $\mathcal{P} \subset \R[X_1,\ldots,X_k]_{\leq d}$ and let  $s = \card({\mathcal{P}})$.
Then
\[ 
\card\left(\bigcup_{\sigma \in \SIGN(\mathcal{P})} \Cc(\RR(\sigma,\R^k))\right)
 \leq
  \sum_{1 \leq j \leq k} \binom{s}{j}
     4^j d (2 d - 1)^{k - 1}. 
\]
\end{proposition}
We now have all the ingredients needed to prove Theorem~\ref{thm:main}.

\subsection{Proof of Theorem~\ref{thm:main}}
\begin{proof}[Proof of Theorem~\ref{thm:main}]
Let $\Phi$ be the $\mathcal{P}$-closed formula defining $S$.
Introducing new variables $Z_1,\ldots,Z_m$, let $\tilde{S} \subset \R^n \times \R^m$ be the $\tilde{\mathcal{P}}$-formula
\[
\Phi \wedge \bigwedge_{1 \leq i \leq m} (Z_i - f_i = 0).
\]
Let $\tilde{f} :\tilde{S} \rightarrow \R^m$ denote the restriction to $\tilde{S}$ of the projection map
$\pi_Z: \R^m \times \R^n \rightarrow  \R^m$ to  the $Z$-coordinates.
Then clearly $S$ is semi-algebraically homeomorphic to $\tilde{S}$, 
$f(S) = \tilde{f}(\tilde{S})$, and $\Reeb(f)$ is semi-algebraically
homeomorphic to $\Reeb(\tilde{f})$. We have the following commutative square where the horizontal arrows are homeomorphisms and the vertical arrows are the quotient maps.
\[
\xymatrix{
S \ar[r]^{\cong}  \ar[d]^\phi&\tilde{S}\ar[d]^{\tilde{\phi}} \\
\Reeb(f) \ar[r]^{\cong} & \Reeb(\tilde{f})
}
\]

Now it follows from Proposition~\ref{prop:parametrized-cc} that there exists a finite set of polynomials $\mathcal{Q} \subset \R[Z_1,\ldots,Z_m]$,
with
\begin{eqnarray}
\label{eqn:Q}
\card(\mathcal{Q}), \max_{Q \in \mathcal{Q}}\deg(Q) &\leq& (sd)^{(n+m)^{O(1)}}
\end{eqnarray}
having the following property:
for each $\sigma \in \SIGN({\mathcal{Q}})$ and each $C \in \Cc(\RR(\sigma,\R^m))$,
there exists
an index set $I_{\sigma,C}$, a finite family of polynomials 
\[
\mathcal{P}_{\sigma,C} \subset \R[X_1,\ldots,X_n,Z_1,\ldots,Z_m],
\]
and 
$\mathcal{P}_{\sigma,C}$ formulas $(\Theta_\alpha(\overline{X},\overline{Z}))_{\alpha \in I_{\sigma,C}}$ 
such that  $\Theta_\alpha(x,z) \Rightarrow z\in C$,  and  
for each $z \in C$, 
and each  connected component $D$ of $\pi_Z^{-1}(C) \cap \tilde{S}$,
there exists a unique $\alpha \in I_{\sigma,C}$ (which does not depend on $z$) with 
$\RR(\Theta_\alpha(\cdot,z)) = \pi_Z^{-1}(z) \cap D$. 

Moreover, the cardinalities of $I_{\sigma,C}$ and $\mathcal{P}_{\sigma,C}$ and the degrees of the polynomials
in $\mathcal{P}_{\sigma,C}$ are all bounded by $(sd)^{(n+m)^{O(1)}}$.

Let $\phi$  (resp.\ $\tilde{\phi}$) be the canonical surjection $\phi: S \rightarrow \Reeb(f) \cong S/\sim$ (resp.\  $\tilde{\phi}: \tilde{S} \rightarrow \Reeb(\tilde{f}) \cong \tilde{S}/\sim$). From Theorem~\ref{thm:definability-of-Reeb} it follows that we can assume that  $\phi$ is a proper semi-algebraic map. For each $i \geq 0$, 
we have the inequality (cf. Proposition~\ref{prop:descent})
\begin{eqnarray}
\label{eqn:descent}
b_i(\Reeb(f)) &\leq& \sum_{p+q=i} b_q(\underbrace{S \times_\phi \cdots \times_\phi S}_{\mbox{$(p+1)$ times}}).
\end{eqnarray}

Now observe  that 
$
\underbrace{\tilde{S} \times_{\tilde{\phi}} \cdots \times_{\tilde{\phi}} \tilde{S}}_{\mbox{$(p+1)$ times}} 
$ 
(and hence $\underbrace{S \times_\phi \cdots \times_\phi S}_{\mbox{$(p+1)$ times}}$)
is semi-algebraically homeomorphic to the semi-algebraic set 
defined by the formula
\begin{eqnarray}
\label{eqn:Theta}
\Theta(\overline{X}^{(0)},\ldots,\overline{X}^{(p)},\overline{Z}) &=& 
\bigvee_{\substack{\sigma \in \SIGN(\mathcal{Q}) \\ C \in \Cc(\RR(\sigma,\R^m)) \\
\alpha \in I_{\sigma,C}}}
\bigwedge_{
 0 \leq j \leq p
 }
\Theta_\alpha(\overline{X}^{(j)},\overline{Z}).
\end{eqnarray}

To see this observe that
\[
((x^{(0)}, z^{(0)}),\ldots, (x^{(p)}, z^{(p)})) \in  \underbrace{\tilde{S} \times_{\tilde{\phi}} \cdots \times_{\tilde{\phi}} \tilde{S}}_{\mbox{$(p+1)$ times}} 
\]
if and only if 
\[
z^{(0)} = \cdots = z^{(p)} = z,
\]
for some $z$, and 
$x^{(0)}, \ldots,x^{(p)}$ belong to the same connected component of $\tilde{f}^{-1}(z)$.

It is easy to verify this last equivalence using the properties of the decomposition given by 
Proposition~\ref{prop:parametrized-cc}.

We now claim that each of the formulas
\[
\Theta(\overline{X}^0,\ldots,\overline{X}^{(p)},\overline{Z}), \; 0 \leq p \leq m,
\]
is a $\tilde{\mathcal{P}}_p$-formula for some finite set $\tilde{\mathcal{P}}_p \subset \R[\overline{X}^0,\ldots,\overline{X}^{(p)},\overline{Z}]$ with 
$\card(\tilde{\mathcal{P}}_p)$ and  the degrees of the polynomials in $\tilde{\mathcal{P}}_p$ being bounded singly exponentially.

In order to prove the claim first observe that
the cardinality of the set
\[
\bigcup_{\sigma \in \SIGN(\mathcal{Q})} \Cc(\RR(\sigma,\R^m))
\] 
is bounded singly exponentially,
once the number of polynomials in $\mathcal{Q}$, and their degrees  are bounded singly exponentially (using Proposition~\ref{7:the:theorem1}).
The fact that the  number of polynomials in $\mathcal{Q}$ and their degrees are bounded singly exponentially follows from
\eqref{eqn:Q}.
Moreover, for similar reasons the cardinalities of the index sets $I_{\sigma, C}$ are also bounded singly exponentially.
The claim now  follows from Eqn. \eqref{eqn:Theta}.

Finally, to prove the theorem we first apply inequality \eqref{eqn:descent} and then apply
Proposition~\ref{prop:GV07} to bound the right hand side of the inequality \eqref{eqn:descent}.
\end{proof}

\begin{remark}
\label{rem:Stein2}
Given the analogy between Reeb spaces and Stein factorization
(cf.\ Remark~\ref{rem:Stein})  it could be interesting to investigate (in the context of algebraic geometry) 
Stein factorization for projective morphisms from the point of view of complexity in analogy with Theorem \ref{thm:main}.  To the best of our knowledge
this has not yet been investigated.
\end{remark}

\section{Conclusion}
\label{sec:conclusion}
In this paper we have proved the realizability of the Reeb space of proper definable maps in an o-minimal structure as
a proper definable quotient. We have exhibited examples where the Reeb spaces of maps can have arbitrarily complicated topology compared to that of the domains of the maps, a sharp contrast with the behavior of Reeb graphs. Nevertheless,
we have proved singly exponential upper bounds on the Betti numbers of the Reeb spaces of proper semi-algebraic maps.


\section{Appendix}
\label{sec:appendix}
This section is devoted to the proof of Proposition~\ref{prop:BElim}.
For the convenience of the reader we first recall the proposition.
\begin{proposition*}
Let $\rR$ be a real closed field.
Let $\mathcal{P} \subset \rR[X_1,\ldots,X_k,Y_1,\ldots,Y_\ell]$ be a finite set of polynomials with degrees bounded by $d$ and with $\card(\mathcal{P}) =s$, and let $S \subset \rR^{k} \times \rR^\ell$ be a bounded $\mathcal{P}$-closed
semi-algebraic set. Then
there exists a finite set of polynomials, $\mathcal{Q} \subset \R[Y_1,\ldots,Y_\ell]$,
with 
\begin{eqnarray*}
\card(\mathcal{Q}) &\leq& (sd)^{(k+\ell)^{O(1)}}, \\
\max_{Q \in \mathcal{Q}} \deg(Q) &\leq& d^{(k+\ell)^{O(1)}},
\end{eqnarray*} 
such that  for  each  
$\sigma \in  \SIGN(\mathcal{Q})$, 
$C \in \Cc(\RR(\sigma,\R^\ell))$,
$\y \in C$, and $D \in \Cc(S_C)$,
$D_\y \in \Cc(S_\y)$.
\end{proposition*}

The rest of this section is devoted to the proof of Proposition~\ref{prop:BElim}. Since the proof is technical and relies on several intermediate results from algorithmic semi-algebraic geometry we first sketch an outline of the proof for the benefit of the reader. 

\subsection{Outline of the proof of Proposition~\ref{prop:BElim}}
\label{subsec:outline-appendix}
The core idea behind the proof comes from \cite{BV06} where a singly exponential upper bound
is proved on the number of homotopy types of the fibers of a semi-algebraic map. 
We explain this idea informally below.

Given $S \subset \R^{k+\ell}$ a $\mathcal{P}$-closed and bounded semi-algebraic set, and $\pi_Y:\R^{k+\ell} \rightarrow \R^\ell$ the projection map onto the last $\ell$-coordinates, we first introduce 
infinitesimal elements $\bar\eps=(\eps_1,\ldots,\eps_s)$ (where $s = \card(\mathcal{P})$), and  define 
a semi-algebraic subset, $S^\star(\bar\eps) \subset \R\la\bar\eps\ra^{k+\ell}$ which  is infinitesimally larger than $S$.
We prove  (see Lemma~\ref{lem:homotopic})  that for each $\y \in \R^\ell$,
the extension of the fiber, $S_\y = S \cap \pi_Y^{-1}(\y) \cap \R^\ell$, to $\R\la\bar\eps\ra$
(see Definition~\ref{def:extension} below for the definition of extension) 
is  a semi-algebraic deformation retract of the fiber $S^\star(\bar\eps)_\y = S^\star(\bar\eps) \cap \pi_Y^{-1}(\y)$ --
and in particular, these two sets fibers are semi-algebraically homotopy equivalent. 
(The field $\R\la\bar\eps\ra$ is a real closed extension of $\R$ consisting of algebraic Puiseux series with coefficients in
$\R$  and is defined below in  Section~\ref{subsec:Puiseux}.) 
It is important to note that in order to obtain the semi-algebraic homotopy equivalence between the 
the extension of the fiber $S_\y$ to $\R\la\bar\eps\ra$  and $S^\star(\bar\eps)_\y$, 
we need the fact that all the coordinates of $\y$ belong to $\R$
(rather than $\R\la\bar\eps\ra \setminus \R$). 
Note that $\R$ is a subfield of $\R\la\bar\eps\ra$, and there is thus a natural inclusion of
$\R^\ell$ into $\R\la\bar\eps\ra^\ell$ as well.

Next, we identify a semi-algebraic subset
$G(\bar\eps) \subset \R\la\bar\eps\ra^{\ell}$  (cf. Definition~\ref{def:criticalvalues})
of codimension at least one in $\R\la\bar\eps\ra^\ell$, 
consisting of the critical values of the projection map restricted to the various algebraic sets 
defined by the polynomials defining $S^\star(\bar\eps)$. 
We prove (Lemma~\ref{lem:discriminant})
that the map $\pi_Y$ restricted to the set 
$\pi_Y^{-1}(\R\la\bar\eps\ra^\ell \setminus G(\bar\eps)) \cap S^\star(\bar\eps)$
is a \emph{locally trivial semi-algebraic
fibration} (see Section~\ref{subsec:fibration} for the definition of semi-algebraic fibrations). 
We also observe that $G(\bar\eps) \cap \R^\ell$ is empty and hence
$\R^\ell \subset \R\la\bar\eps\ra^\ell \setminus G(\bar\eps)$.

Moreover,  for each 
$C \in \Cc(\R\la\bar\eps\ra^\ell \setminus G(\bar\eps))$, 
$C \cap \R^\ell$ is a semi-algebraic subset of $\R^\ell$. If the projection
map $\pi_Y$ restricted to $\pi_Y^{-1}(C) \cap S^\star(\bar\eps)$ is a \emph{trivial}  semi-algebraic fibration (not just a \emph{locally trivial} one), 
then we would be done at this point, because we could take $\mathcal{Q}\subset \R[Y_1,\ldots,Y_\ell]$ to
be a finite set of polynomials such that the semi-algebraic sets 
\[
C \cap \R^\ell, C \in \Cc(\R\la\bar\eps\ra^\ell \setminus  G(\bar\eps))
\] 
are all $\mathcal{Q}$-semi-algebraic sets. (It follows from certain preliminary results (Propositions~\ref{prop:cc} and \ref{prop:qe}) that such a set
$\mathcal{Q}$ exists satisfying the singly exponential estimates in Proposition~\ref{prop:BElim}.)

In order to arrive at the desired situation (i.e to have a trivial semi-algebraic fibration instead of just a locally trivial one) 
we need an additional step of covering the complement of $G(\bar\eps)$ using semi-algebraically contractible sets noting that a locally trivial semi-algebraic fibration over a contractible semi-algebraic set is automatically trivial (see Proposition~\ref{prop:fibration1}).
We use a result from \cite{BPRbook2} which gives an algorithm for constructing a cover of a given semi-algebraic set by contractible ones (see Proposition~\ref{prop:covering}) whose complexity is bounded singly exponentially.
This introduces some additional technicalities since we need to consider an infinitesimal thickening of
$G(\bar\eps)$ before taking a cover by contractible semi-algebraic sets  
of its complement and this is explained in the proof.  

One important technical detail is that we need to keep track of the dependence of each polynomial
appearing in the definition of the cover on the various infinitesimals introduced in the course of the construction. It is important that each such polynomial depends on at most polynomially many of these
infinitesimals and the degrees of the polynomials in the infinitesimals are bounded singly exponentially
(cf. Parts~\eqref{itemlabel:complexity-of-H} and \eqref{itemlabel:dependence-of-H} of Lemma~\ref{lem:H}).
This is crucial for  ensuring  the singly exponential bounds on the set of polynomials $\mathcal{Q}$ that is produced at the end.

Before proving Proposition~\ref{prop:BElim} we need a few preliminary results.

\subsection{Preliminaries}
\label{subsec:preliminaries}
\subsubsection{Real closed extensions and Puiseux series}
\label{subsec:Puiseux}
We will need some
properties of Puiseux series with coefficients in a real closed field. We
refer the reader to \cite{BPRbook2} for further details.

\begin{notation}
\label{not:Puiseux}
  For $\rR$ a real closed field we denote by $\rR\la \eps\ra$ 
  the real closed field of algebraic Puiseux series in $\eps$
  with coefficients in $\rR$. We use the notation $\rR\langle \eps_{1},
  \ldots, \eps_{m} \rangle$ to denote the real closed field $\R\langle \eps_{1}\rangle \langle \eps_{2}\rangle
  \cdots \langle \eps_{m}\rangle$. Note that in the unique
  ordering of the field $\rR \langle \eps_{1}, \ldots, \eps_{m}
  \rangle$, $0< \eps_{m} \ll \eps_{m-1} \ll \cdots \ll \eps_{1} \ll 1$.
\end{notation}

\begin{notation}
\label{not:lim}
  For elements $x \in \rR\langle \eps\rangle$ which are bounded
  over $\rR$ we denote by $\lim_{\eps}  x$ to be the image in $\rR$ under the
  usual map that sets $\eps$ to $0$ in the Puiseux series $x$.
  If $\bar\eps =(\eps_1,\ldots,\eps_m)$, and $x \in \rR \langle \bar\eps\rangle$ bounded
  over $\rR$, then we will denote
  \[
  \lim_{\bar\eps} x = \lim_{\eps_1}( \cdots (\lim_{\eps_m} x)).
  \]
\end{notation}

\begin{definition}
\label{def:extension}
  If $\rR'$ is a real closed extension of a real closed field $\rR$, and $S
  \subset \rR^{k}$ is a semi-algebraic set defined by a first-order formula
  with coefficients in $\rR$, then we will denote by $\E(S, \rR') \subset \rR'^{k}$ the semi-algebraic subset of $\rR'^{k}$ defined by
  the same formula. It is well-known that $\E(S, \rR')$ does
  not depend on the choice of the formula defining $S$ \cite{BPRbook2}.
\end{definition}

For $\phi$ a $\mathcal{P}$-closed formula, 
with $\mathcal{P} = \{P_1,\ldots,P_s\}$, 
we will denote by $\phi^\star(\bar\delta)$ the formula obtained from $\phi$
by replacing any occurrence of 
  $P_{i} \ge 0$  with $P_{i} \ge -\delta_i$, and
  any occurrence of  $P_{i} \le 0$  with $P_{i} \le \delta_i$, for each $i, 1 \leq i \leq s$.

\begin{lemma}
  \label{16:lem:star} 
  Let $\R$ be a real closed field and $R \in \R$ with $R > 0$.
  The semi-algebraic set 
  $\E(\RR(\phi,\R^k) \cap [-R,R]^k,
  \R \la \bar\delta \ra)
  $ is semi-algebraically
  homotopy equivalent to $\RR(\phi^\star(\bar\delta), \R\la\bar\delta\ra^k) \cap 
  [-R,R]^k$.
\end{lemma}

\subsubsection{Covering by semi-algebraic contractible sets}
\label{subsec:cover}
As explained in the outline we will need an auxiliary result on covering of semi-algebraic sets by
contractible ones with singly exponential complexity. 
In order to state this result we first need the notion of being in \emph{strong general position}.

\begin{definition}
\label{def:general-position}
Let $\rR$ be a real closed field and 
$\mathcal{P}  \subset \rR [X_{1} , \ldots ,X_{k}
]$ be a  finite set. We say that the family $\mathcal{P}$ is in
\emph{$\ell$-general position}, if no more than 
$\ell$ polynomials belonging to $\mathcal{P}$ have a zero in $\rR^{k}$.
The family $\mathcal{P}$ is in \emph{strong $\ell$-general
position}  if moreover any $\ell$
polynomials belonging to $\mathcal{P}$ have at most a
finite number of zeros in $\rR^{k}$.
\end{definition}

We will use the following proposition which follows from the correctness and complexity analysis of 
Algorithm 16.14 (Covering by Contractible Sets)  in  \cite{BPRbook2}.

\begin{proposition}
\label{prop:covering}
Let $\rR$ be a real closed field and $\D \subset \rR$ and ordered domain.
There exists an algorithm that takes as input
\begin{enumerate}[(1)]
\item  
a finite set of $s$ polynomials $\mathcal{G}= \{G_1,\ldots,G_t \}  \subset \D[\bar\eps,\bar\delta][Y_1,\ldots,Y_\ell]$,
where $\bar\eps = (\eps_1,\ldots,\eps_s)$, $\bar\delta = (\delta_1,\ldots,\delta_t)$,
in strong $\ell$-general position on $\rR\la\bar\eps,\bar\delta\ra^\ell$,
and such that each polynomial in $\mathcal{G}$ depends on at most $k+\ell$ of the $\eps_i$'s and at most one of the $\delta_i$'s;
\item a  $\mathcal{G}$-closed formula $\psi$;
\item $R > 0$, $R \in \mathrm{D}$;
\end{enumerate}
and outputs 
\begin{enumerate}[(a)]
\item a finite set of polynomials $\mathcal{H} \subset \D[\bar\eps,\bar\delta,\bar\zeta][Y_1,\ldots,Y_\ell]$,
where $\bar\zeta = (\zeta_1,\ldots,\zeta_{2 \card(\mathcal{H})})$;
\item
a tuple  of $\mathcal{H}$-formulas $(\theta_\alpha)_{\alpha \in I}$ such that
each  $\RR(\theta_\alpha,\rR\la\bar\eps,\bar\delta,\bar\zeta\ra,^\ell), \alpha \in I$ is a closed
semi-algebraically contractible set, and
\item 
\[
\bigcup_{\alpha \in I} 
\RR(\theta_\alpha,\rR\la\bar\eps,\bar\delta,\bar\zeta\ra^\ell) = \RR(\psi,\rR\la\bar\eps,\bar\delta,\bar\zeta\ra^\ell)
\cap [-R,R]^{\ell}.
\]
\end{enumerate}
Moreover,
\begin{eqnarray*}
\card(I), \card(\mathcal{H}) &\leq& (t D)^{\ell^{O(1)}}, \\ 
\deg_{\bar{Y}}(H), \deg_{\bar{\eps}}(H), \deg_{\bar{\delta}}(H),\deg_{\bar{\zeta}}(H) &\leq&
D^{\ell^{O(1)}},
\end{eqnarray*}
where $D = \max_{1 \leq i \leq t} \deg(G_i)$.

Moreover, each polynomial appearing in $\mathcal{H}$ depends on at most
$(k+\ell) (\ell+1)^2$ of $\eps_i$'s, at most $(\ell+1)^2$ of the $\delta_i$'s and on at most one of the $\zeta_i$'s. 
\end{proposition}

\begin{remark}
\label{rem:local}
Note that the last claim in Proposition~\ref{prop:covering}, namely that each polynomial appearing in any of the formulas $\theta_\alpha$ depends on at most
$(k+\ell) (\ell+1)^2$ of $\eps_i$'s, at most $(\ell+1)^2$ of the $\delta_i$'s and on at most one of the $\zeta_i$'s,
does not appear explicitly in \cite{BPRbook2}, but is evident on a close examination of the algorithm. It is also 
reflected in the fact that the combinatorial part (i.e. the part depending on $\card(\mathcal{G})$) 
of the complexity of Algorithm 16.14 in \cite{BPRbook2} is bounded by $\card(\mathcal{G})^{(\ell+1)^2}$. This is
because the Algorithm 16.14 in \cite{BPRbook2} has a ``local property'', namely that all computations involve
at most a small number  (in this case $(\ell+1)^2$) polynomials in the input at a time.
\end{remark}

\subsubsection{Effective quantifier elimination}
\label{subsec:qe}
We will also need the following effective bound on the complexity of eliminating one block of 
existential quantifiers in the theory of real closed fields. It is a direct consequence of Theorem 14.16 in \cite{BPRbook2}.
\begin{proposition}
\label{prop:qe}
Let $\rR$ be a real closed field, $\D \subset \rR$ an ordered domain,  
$\mathcal{Q} \subset \D[\eps_1,\ldots,\eps_m][X_1,\ldots,X_k,Y_1,\ldots,Y_\ell]$ a finite set of polynomials, and $\psi$ a $\mathcal{Q}$-formula. 
Suppose 
that the degrees of the polynomials in $\mathcal{Q}$ in $\eps_h, X_i, Y_j$ for $1 \leq h \leq m, 1\leq i \leq k, 1\leq j \leq \ell$  are all bounded by $D$.

Then, there exists a finite set of polynomials
$\mathcal{G} \subset \D[\eps_1,\ldots,\eps_m][Y_1,\ldots,Y_\ell]$, and a $\mathcal{G}$-formula $\theta$ satisfying the
following conditions:
\begin{enumerate}[(a)]
\item
\label{itemlabel:prop:qe:a}
\[
\card(\mathcal{G}) = (\card(\mathcal{Q}) D)^{O((k+1)(\ell+1))};
\]
\item
\label{itemlabel:prop:qe:b}
\[
\max_{G \in \mathcal{G}, 1\leq i \leq m, 1 \leq j \leq \ell} (\deg_{\eps_i}(G), \deg_{Y_j}(G)) = D^{O(k)};
\]
\item
\label{itemlabel:prop:qe:c}
\[
\RR(\theta,\rR\la\bar\eps\ra^\ell) = \pi_Y(\RR(\psi,\rR\la\bar\eps\ra^{k+\ell})).
\]
\end{enumerate}
\end{proposition}

\subsubsection{Trivial and locally trivial semi-algebraic fibrations}
\label{subsec:fibration}
As explained in the outline we will also need the notion of trivial and locally trivial semi-algebraic fibrations.

\begin{definition}
Let $\rR$ be a real closed field. 
Let $p: E \rightarrow B$ be a semi-algebraic map.

We say that $p$ is a \emph {trivial semi-algebraic  fibration}
if there exists $b \in B$ and a semi-algebraic homeomorphism $h: E \rightarrow p^{-1}(b) \times B$, such that
the following diagram commutes:
\[
\xymatrix{
E \ar[rr]^{h} \ar[rd]^{p} && p^{-1}(b) \times B \ar[ld]^{\pi} \\
&B&
}
\]
(here $\pi$ is the projection to the second factor).

We say that the map $p$ is a \emph{locally trivial semi-algebraic fibration} if there exists a finite
covering $(U_i)_{i \in I}$ of $B$ by semi-algebraic subsets which are open in $B$, such that
for each $i \in I$, $p|_{p^{-1}(B_i)}$ is a semi-algebraic trivial fibration.
\end{definition}

\begin{definition}[Pull-back of a semi-algebraic locally trivial fibration under a semi-algebraic map]
Let $p:E \rightarrow B$ be a locally trivial semi-algebraic fibration, and $f:B' \rightarrow B$ be a 
continuous semi-algebraic map. Then the map $f^*p: E \times_B B' \rightarrow B',  (e,b') \mapsto b'$
is again a semi-algebraic locally trivial fibration (called the \emph{pull-back of $p$ under $f$}).
\end{definition}

The main property of semi-algebraic fibrations that we will need is the following.
\begin{proposition}
\label{prop:fibration1}
Let $p: E \rightarrow B$ be a  locally trivial semi-algebraic fibration. Suppose $B'$ is a semi-algebraically
contractible subset of $B$, and $B''$ any semi-algebraic subset of $B'$. Then $p|_{p^{-1}(B'')}$ is a  trivial semi-algebraic fibration.  
\end{proposition}

Proposition~\ref{prop:fibration1} is a corollary of the following more general proposition.

\begin{proposition}
\label{prop:fibration2}
Let $p: E \rightarrow B$ be a  locally trivial semi-algebraic fibration, and let $f,g: B' \rightarrow B$ be two continuous semi-algebraic maps which are semi-algebraically homotopic.
The, $f^*p$ and $g^*p$ are isomorphic as locally trivial semi-algebraic fibrations.
\end{proposition}

\begin{proof}
The corresponding statement in the  category of topological spaces and maps appear in 
\cite[Theorem 4.6.4]{Aguilar} with the 
extra assumption that $B'$ is paracompact. 
We do not give a complete proof of Proposition~\ref{prop:fibration2} here but just indicate the modifications needed in the proof of Theorem 4.6.4 in \cite{Aguilar}
to carry it over to the semi-algebraic setting over an arbitrary real closed field $\R$.

The key result used in the proof of Theorem 4.6.4 in \cite{Aguilar} is Proposition 4.6.3 in \cite{Aguilar} whose statement is reproduced below using the notation of the current paper.

``Let $p:E \rightarrow B \times [0,1]$ be a locally trivial fibration, where $B$ is a paracompact space. Let 
$r:B \times [0,1] \rightarrow B \times [0,1]$ be the retraction defined by $r(b,t) =  (b,1)$. Then, 
$p$ is isomorphic to $r^*p$.''

We will need the following semi-algebraic version of the above statement, namely:

\begin{claim}
\label{claim:proof:prop:fibration2:1}
Let $p:E \rightarrow B \times [0,1]$ be a locally trivial semi-algebraic fibration. Let 
$r:B' \times [0,1] \rightarrow B \times [0,1]$ be the retraction defined by $r(b,t) =  (b,1)$. Then, 
$p$ is isomorphic to $r^*p$.
\end{claim}

All the modifications needed are in the proof of the above claim.
In the proof of Proposition 4.6.3 in \cite{Aguilar}
replace the open cover $(U_\alpha)_{\alpha  \in \Lambda}$ of $B$ (which is assumed to be only locally finite using the  paracompactness of $B$)
by a finite open cover of $(U_i)_{i \in I }$ of the semi-algebraic set  
$B$ such that $p$ restricted to each $U_\alpha \times [0,1]$ is a
trivial semi-algebraic fibration.  Such a finite cover exists if 
$p$ is a locally trivial semi-algebraic fibration. Next, replace the partition of unity $\{\eta_\alpha\}_{\alpha \in \Lambda}$ subordinate to the partition $(U_\alpha)_{\alpha  \in \Lambda}$ by a 
a semi-algebraic partition of unity $\{\eta_i\}_{i \in I}$ subordinate to the (finite) cover
$(U_i)_{i \in I }$ which is known to exist (see \cite[Lemma 12.7.3]{BCR}). These two modifications to the proof of 
Proposition 4.6.3 in \cite{Aguilar} produce a proof of Claim~\ref{claim:proof:prop:fibration2:1}. 
Proposition~\ref{prop:fibration2} now follows from  Claim~\ref{claim:proof:prop:fibration2:1}
in exactly the same way as Theorem 4.6.4 is deduced from Proposition 4.6.3 in \cite{Aguilar}.
\end{proof}

\begin{proof}[Proof of Proposition~\ref{prop:fibration1}]
First observe that it is obvious that the restriction of a trivial semi-algebraic fibration $p':E'\rightarrow B'$ to
to a semi-algebraic subset $B'' \subset B'$ is again a trivial semi-algebraic fibration. 
Let $p':E' \rightarrow B'$ be the restriction of the locally trivial semi-algebraic fibration $p$ to $B'$.
Now since $B'$ is assumed to semi-algebraically contractible, there exists a semi-algebraic homotopy
between $f= \mathrm{id}_{B'}$ and a constant map $g:B' \rightarrow B'$. We obtain that
$p' \cong f^*p' \cong g^* p'$, where the first isomorphism is due to the fact that $f$ is the identity map, and 
the second isomorphism is a consequence of Proposition~\ref{prop:fibration2} and the fact that $f$ and $g$ are homotopic.
Now since $g^*p'$ is the pull-back of $p'$ under a constant map, it is clearly a trivial
semi-algebraic fibration. Finally, the observation at the beginning of the proof implies that $p'$ restricted to
$B'' \subset B'$ is also trivial.
\end{proof}

\subsection{Proof of Proposition~\ref{prop:BElim}}
\label{subsec:proof}
 
For the rest of this section we fix a real closed field $\R$, a finite set of polynomials
\[
\mathcal{P}  =\{P_1,\ldots,P_s\} 
\subset \rR[X_1,\ldots,X_k,Y_1,\ldots,Y_\ell],
\]
as well as a $\mathcal{P}$-closed formula
$\phi$. We denote $S = \RR(\phi,\rR^{k+\ell})$.

\begin{notation}
\label{not:start}
For $\bar{\eps} = (\eps_1,\ldots,\eps_{s})$, we denote by
$\phi^{\star,c} (\bar{\eps})$, the $\mathcal{P}^\star(\bar{\eps})$-closed formula
obtained by replacing each occurrence of $P_i \geq 0$ in $\phi$ by $P_i + \eps_i \geq 0$ 
(resp. $P_i \leq 0$ in $\phi$ by $P_i - \eps_i \leq 0$) for $1 \leq i \leq s$,
where 
\[
\mathcal{P}^\star(\bar{\eps}) = \bigcup_{1 \leq i \leq s} \{ P_i +\eps_i, P_i - \eps_i \}.
\]
\end{notation}

Observe that
\begin{equation}
\label{eqn:S-star}
S^{\star,c} (\bar{\eps}):= \RR(\phi^{\star,c} (\bar{\eps}), \rR\la\bar{\eps}\ra^{k+\ell})
\end{equation}
is a $\mathcal{P}^\star(\bar{\eps})$-closed semi-algebraic set.

\begin{remark}
\label{rem:Puiseux}
We will use a few results whose proofs can be found in \cite{BV06}. These results are stated in \cite{BV06}
not using the language of Puiseux extensions but rather with the $\eps_i$'s occurring in the statement 
as small enough positive elements of the ground field $\rR$. 
However, the statements in \cite{BV06}  imply the corresponding statements of the
current paper in terms of Puiseux series 
by an application of a standard argument using the Tarski-Seidenberg transfer principle.

More precisely, the following fact  which  is a consequence of the Tarski-Seidenberg transfer principle
will be used repeatedly.

If $\phi(T)$ is a first order formula with constants in a real closed field $\rR$,  
then the first order sentence 
\[
(\exists T_0) (T_0 > 0) \wedge (\forall T) ((0 < T) \wedge (T < T_0)) \Rightarrow \phi(T)
\]   
is true over $\rR$ if and only if
the sentence 
$\phi(\eps)$ is true over $\rR\la\eps\ra$. 

In order to see this observe that $\RR(\phi(T),\rR)$ is a semi-algebraic subset of $\rR$ which contains an interval
$(0,t_0)$ for some $t_0 > 0$. 
Then by the Tarski-Seidenberg transfer principle, $(0,t_0) \subset \RR(\phi(T),\rR\la\eps\ra)$ as well.
Now, since $\eps$ is positive and smaller than all positive elements of $\rR$ in the unique ordering
of the real closed field $\rR\la\eps\ra$, $0 < \eps < t_0$, 
and so $\eps \in \RR(\phi(T),\rR\la\eps\ra)$.
Hence, $\phi(\eps)$ is true over $\R\la\eps\ra$. 
\end{remark}

\begin{lemma}
\label{lem:non-singular}
For each $\mathcal{Q} \subset \mathcal{P}^\star(\bar\eps)$,
$\ZZ(\mathcal{Q},\rR\la\bar\eps\ra^{k+\ell})$ is either empty or 
is a non-singular $(k+\ell-\card(\mathcal{Q}))$-dimensional 
real variety
such that at every point 
\[
(\x,\y) = (x_1,\ldots,x_k,y_1,\ldots,y_\ell) \in \ZZ(\mathcal{Q},\rR\la\bar\eps\ra^{k+\ell}),
\] 
the $(\card(\mathcal{Q}) \times (k+\ell))$-Jacobian matrix,
\[
\left( \frac{\partial P}{\partial X_i} , \frac{\partial P}{\partial Y_j}
\right)_{\substack{P \in \mathcal{Q} \\
1\leq i \leq k, 1\leq j \leq \ell}}
\]
has the maximal possible rank.
\end{lemma}

\begin{proof}
See proof of Lemma 3.3 in \cite{BV06} and use Remark~\ref{rem:Puiseux}.
\end{proof}

\begin{lemma}
\label{lem:non-singular2}
For each $\mathcal{Q} \subset \mathcal{P}^\star(\bar\eps)$,
and $\y \in \rR^\ell$,
$\ZZ(\mathcal{Q}(\cdot,\y),\rR\la\bar\eps\ra^{k})$ is either empty or 
is a non-singular $(k-\card(\mathcal{Q}))$-dimensional 
real variety
such that at every point 
\[
(\x,\y) = (x_1,\ldots,x_k,y_1,\ldots,y_\ell) \in \ZZ(\mathcal{Q},\rR\la\bar\eps\ra^{k+\ell}),
\] 
the $(\card(\mathcal{Q}) \times k)$-Jacobian matrix,
\[
\left( \frac{\partial P}{\partial X_i}
\right)_{\substack{P \in \mathcal{Q} \\
1\leq i \leq k}}
\]
has the maximal possible rank.
\end{lemma}

\begin{proof}
Similar to proof of Lemma~\ref{lem:non-singular}.
\end{proof}

\begin{notation}[Critical points and critical values]
\label{not:crit}
For $\mathcal{Q} \subset \mathcal{P}^\star(\bar\eps)$, we denote by $\Crit(\mathcal{Q})$
the subset of $\ZZ(\mathcal{Q},\rR\la\bar\eps\ra^{k+\ell})$ at which the Jacobian matrix,
\[
\left( \frac{\partial P}{\partial X_i} \right)_{\substack{P \in \mathcal{Q}, \\
1 \leq i \leq k}}
\]
is not of the maximal possible rank.
We denote $\crit(\mathcal{Q}) = \pi_Y(\Crit(\mathcal{Q}))$.
\end{notation}

\begin{definition}
\label{def:criticalvalues}
Let 
\[
G_1(\bar\eps) =  \bigcup_{\substack{\mathcal{Q} \subset \mathcal{P}^\star(\bar\eps)\\
0< \card(\mathcal{Q}) \leq k}}  \crit(\mathcal{Q}),
\]
and
\[
G_2(\bar\eps) =  \bigcup_{\substack{\mathcal{Q} \subset \mathcal{P}^\star(\bar\eps) \\k < \card(\mathcal{Q}) \leq k +\ell}}  \pi_Y(\ZZ(\mathcal{Q},\rR\la\bar\eps\ra^{k+\ell})),
\]
and
\[
G(\bar\eps) = G_1(\bar\eps) \cup G_2(\bar\eps).
\]
\end{definition}

\begin{lemma}
\label{lem:discriminant}
Let 
\[
\widehat S:=  \pi_Y^{-1}(\R\la\bar\eps\ra^\ell \setminus G(\bar\eps))  \cap S^{\star,c}(\bar\eps).
\] 

Then, the  map $\pi_Y|_{\widehat{S}}$ is a locally trivial  semi-algebraic fibration.
\end{lemma}

\begin{proof}
See proof of Lemma 3.8 in \cite{BV06} and use Remark~\ref{rem:Puiseux}.
\end{proof}

\begin{lemma}
\label{lem:G}
There exists a finite set of polynomials $\mathcal{G} \{G_1,\ldots,G_t \} \subset \R[\bar\eps][Y_1,\ldots,Y_\ell]$, and 
a $\mathcal{G}$-formula $\psi$ satisfying the following.
\begin{enumerate}[(a)]
\item
\[
G(\bar\eps) = \RR(\psi,\rR\la\bar\eps\ra^\ell),
\]
\item
\label{itemlabel:card-of-G}
\[
\card(\mathcal{G}) \leq (k s d)^{O(k\ell)},
\]
\item
\label{itemlabel:deg-of-G}
\[
\max_{1 \leq i \leq t} (\deg_{\bar{Y}}(G_i), \deg_{\bar\eps}\deg(G_i)) \leq  d^{O(k\ell)}.
\]
\item
\label{itemlabel:dependence-of-G}
Moreover, each  polynomial in $\mathcal{G}$ 
depends on at most $(k+\ell)$ of the $\eps_i$'s.
\end{enumerate}
\end{lemma}

\begin{proof}
For each subset $\mathcal{Q} \subset \mathcal{P}^\star(\bar\eps)$ with $0< \card(\mathcal{Q}) \leq k$,
use Proposition~\ref{prop:qe} to obtain a set of polynomials
$\mathcal{G}_{1,\mathcal{Q}}  \subset \rR[\bar\eps][Y_1,\ldots,Y_\ell]$  such that
$\crit(\mathcal{Q})$ is a $\mathcal{G}_{1,\mathcal{Q}}$-semi-algebraic set defined by  a $\mathcal{G}_{1,\mathcal{Q}}$-formula 
formula $\psi_{1,\mathcal{Q}}$.

For each subset $\mathcal{Q} \subset \mathcal{P}^\star(\bar\eps)$ with $k< \card(\mathcal{Q}) \leq k+\ell$,
use Proposition~\ref{prop:qe} to obtain a set of polynomials
$\mathcal{G}_{2,\mathcal{Q}}  \subset \rR[\bar\eps][Y_1,\ldots,Y_\ell]$  such that
$\pi_Y(\ZZ(\mathcal{Q},\rR\la\bar\eps
\ra^\ell))$ is a $\mathcal{G}_{2,\mathcal{Q}}$-semi-algebraic set defined by  a $\mathcal{G}_{2,\mathcal{Q}}$-formula 
formula $\psi_{2,\mathcal{Q}}$.

Let 
\[
\mathcal{G} = \{G_1,\ldots,G_t\} =  \bigcup_{\substack{\mathcal{Q} \subset \mathcal{P}^\star(\bar\eps)\\
0< \card(\mathcal{Q}) \leq k}}   \mathcal{G}_{1,\mathcal{Q}}\cup \bigcup_{\substack{\mathcal{Q} \subset \mathcal{P}^\star(\bar\eps) \\k < \card(\mathcal{Q}) \leq k +\ell}}  \mathcal{G}_{2,\mathcal{Q}},
\]

and 
\[
\psi=   \bigvee_{\substack{\mathcal{Q} \subset \mathcal{P}^\star(\bar\eps)\\
0< \card(\mathcal{Q}) \leq k}}   \psi_{1,\mathcal{Q}}\vee \bigvee_{\substack{\mathcal{Q} \subset \mathcal{P}^\star(\bar\eps) \\k < \card(\mathcal{Q}) \leq k +\ell}}  \psi_{2,\mathcal{Q}}.
\]

Then, 
\begin{eqnarray*}
\RR(\psi,\R\la\bar\eps\ra^\ell) &=&  \bigcup_{\substack{\mathcal{Q} \subset \mathcal{P}^\star(\bar\eps)\\
0< \card(\mathcal{Q}) \leq k}}   \RR(\psi_{1,\mathcal{Q}},\R\la\bar\eps\ra^\ell)  \cup \bigcup_{\substack{\mathcal{Q} \subset \mathcal{P}^\star(\bar\eps) \\k < \card(\mathcal{Q}) \leq k +\ell}}  
\RR(\psi_{2,\mathcal{Q}},\R\la\bar\eps\ra^\ell) \\
&=& 
\bigcup_{\substack{\mathcal{Q} \subset \mathcal{P}^\star(\bar\eps)\\
0< \card(\mathcal{Q}) \leq k}}   \pi_Y(\Crit(\mathcal{Q}))  \cup \bigcup_{\substack{\mathcal{Q} \subset \mathcal{P}^\star(\bar\eps) \\k < \card(\mathcal{Q}) \leq k +\ell}}  
\pi_Y(\ZZ(\mathcal{Q},\rR\la\bar\eps
\ra^\ell)) \\
&=&
G_1(\bar\eps) \cup G_2(\bar\eps) \\
&=& G(\bar\eps).
\end{eqnarray*}

Observe that each of the sets $\mathcal{G}_{i,\mathcal{Q}}$, where $i=1,2$ and $\mathcal{Q} \subset \mathcal{P}^*(\bar\eps)$, 
depends on at most $(k+\ell)$ of the $\eps_i$'s.
Moreover, using  Parts~\eqref{itemlabel:prop:qe:b} and \eqref{itemlabel:prop:qe:c} of Proposition~\ref{prop:qe} we can assume that 
\begin{eqnarray*}
\card(\mathcal{G}) &\leq& (k s d)^{O(k\ell)},\\
\max_{G \in \mathcal{G}} (\deg_Y(G), \deg_{\bar\eps}\deg(G)) &\leq&  d^{O(k\ell)}
\end{eqnarray*} 
(where 
$s = \card(\mathcal{P})$ and $d = \max_{P \in \mathcal{P}} \deg(P)$).
\end{proof}

Without loss of generality we can assume that 
\[
\psi = \bigvee_{\sigma \in \Sigma} \psi_\sigma,
\]
for some subset $\Sigma \subset \SIGN(\mathcal{G})$, where for $\sigma \in \Sigma$,
$\psi_\sigma$ is the formula defined by
\[
\psi_\sigma = \bigwedge_{1 \leq i \leq t} (G_i \sigma_i 0),
\]
where $\sigma_i$ equals $>,<,0$ according to  $\sigma(G_i) = 1,-1,0$, respectively.

Denote by
$\psi^{\star,o} (\bar\eps,\bar\delta)$, the $\mathcal{G}^\star(\bar\eps,\bar\delta)$-open formula
obtained by replacing each occurrence of $G_i = 0$ in $\psi$ by $(G_i - \delta_i < 0) \wedge (G_i + \delta_i > 0)  $,
$G_i >  0$ in $\psi$ by $G_i + \delta_i > 0$,
and each occurrence of $G_i < 0$ by $G_i - \delta_i  < 0$,
for $1 \leq i \leq t$,
where 
\[
\mathcal{G}^\star(\bar\eps,\bar\delta) = \bigcup_{1 \leq i \leq t} \{ G_i + \delta_i, G_i - \delta_i \}.
\]

Let
\[
G^{+,o}(\bar\eps,\bar\delta) := \RR(\psi^{\star,o} (\bar\eps,\bar\delta), \R\la\bar\eps,\bar\delta\ra^\ell).
\] 

Then,
\[
G^{+,o}(\bar\eps,\bar\delta)^c := \R\la\bar\eps,\bar\delta\ra^\ell \setminus G^{+,o}(\bar\eps,\bar\delta)
\] 
is 
a $\mathcal{G}^\star(\bar\eps,\bar\delta)$-closed semi-algebraic set. 

\begin{lemma}
\label{lem:general-position}
The set $\mathcal{G}^\star(\bar\eps,\bar\delta)$ is
in strong $\ell$-general position over $\rR\la\bar\eps,\bar\delta\ra^\ell$.
\end{lemma}

\begin{proof}
Similar to proof of Lemma~\ref{lem:non-singular}.
\end{proof}

\begin{lemma}
\label{lem:H}
For each $R > 0, R \in \rR$,
there exists a finite set of polynomials $\mathcal{H} \subset \rR[\bar\eps,\bar\delta,\bar\zeta][Y_1,\ldots,Y_\ell]$,
where
$\bar\zeta = (\zeta_1,\ldots,\zeta_{2\card(\mathcal{H})})$,
and a tuple of $(C_\alpha)_{I \in \alpha}$ of $\mathcal{H}$-closed semi-algebraic sets satisfying the
following.
\begin{enumerate}[(a)]
\item
\label{itemlabel:contractibility-of-H}
Each $C_\alpha$ is a
$\mathcal{H}$-semi-algebraic set which is moreover semi-algebraically contractible.
\item
\label{itemlabel:coverproperty-of-H}
\[
\bigcup_{\alpha \in I} C_\alpha =  G^{+,o}(\bar\eps,\bar\delta)^c\cap [-R,R]^\ell.
\]
\item
\label{itemlabel:complexity-of-H}
\begin{eqnarray*}
\card(I), \card(\mathcal{H}) &\leq& (s d)^{k^{O(\ell)}}, \\
\deg_{\bar{Y}}(H), \deg_{\bar{\eps}}(H), \deg_{\bar{\delta}}(H), \deg_{\bar{\zeta}}(H) &\leq& d^{k^{O(\ell)}}.
\end{eqnarray*}

\item
\label{itemlabel:dependence-of-H}
Moreover, each polynomial in $\mathcal{H}$ depends on at most
$(k+\ell)(\ell+1)^2$ of $\eps_i$'s, at most $(\ell+1)^2$ of the $\delta_i$'s and on at most one of the $\zeta_i$'s. 
\end{enumerate}
\end{lemma}

\begin{proof}
Use  Lemma~\ref{lem:G} and Proposition~\ref{prop:covering}.
\end{proof}

Since the semi-algebraic set $S = \RR(\phi,\rR^{k+\ell})$ is assumed to be bounded over $\rR$, 
there exists $R \in \rR, R > 0$, such that
$S \subset [-R,R]^{k+\ell}$. We fix $R > 0$, such that $S \subset [-R,R]^{k+\ell}$
in what follows.
We also fix the cover $(C_\alpha)_{\alpha \in I}$ and the finite family of polynomials $\mathcal{H} \subset \rR[\bar\eps,\bar\delta,\bar\zeta][Y_1,\ldots,Y_\ell]$
given by Lemma~\ref{lem:H}.

\begin{lemma}
\label{lem:trivial-fibration}
For each $\sigma \in \SIGN(\mathcal{H})$, such that
$\RR(\sigma,\rR\la\bar\eps,\bar\delta,\bar\zeta\ra^\ell) \subset C_\alpha$ for some
$\alpha \in I$,  
let
\[
S_\sigma:=  \pi_Y^{-1}(\RR(\sigma,\rR\la\bar\eps,\bar\delta,\bar\zeta\ra^\ell) )  \cap \E(S^\star(\bar\eps), \rR\la\bar\eps,\bar\delta,\bar\zeta\ra).
\] 
Then, the  map $\pi_Y|_{S_\sigma}$ is a trivial  semi-algebraic fibration.
\end{lemma}

\begin{proof}
Follows from Lemma~\ref{lem:discriminant} and Proposition~\ref{prop:fibration1} noting that the 
semi-algebraic set $C_\alpha$ is semi-algebraically contractible using Part \eqref{itemlabel:contractibility-of-H}
of Lemma~\ref{lem:H}.
\end{proof}

Using the same notation as above we have:

\begin{lemma}
\label{lem:Q}
\begin{enumerate}[(a)]
\item
\label{itemlabel:lem:Q:1}
\[
\rR^\ell  \cap [-R,R]^\ell   \subset \bigcup_{\substack{\sigma \in \SIGN(\mathcal{H}), 
\\
\RR(\sigma) \subset \bigcup_{\alpha \in I} C_\alpha}}
\RR(\sigma,\rR\la\bar\eps,\bar\delta,\bar\zeta\ra^\ell).
\]
\item
\label{itemlabel:lem:Q:2}
There exists  a finite set of polynomials, $\mathcal{Q} \subset \rR[Y_1,\ldots,Y_\ell]$,
with degrees and cardinality bounded by $(sd)^{(k+\ell)^{O(1)}}$, 
such that for each 
$\sigma \in  \SIGN(\mathcal{Q})$,
there exists $\alpha \in I$, such that $\RR(\sigma,\rR^\ell) \subset C_\alpha$. 
\end{enumerate}
\end{lemma}

\begin{proof}

\begin{proof}[Proof of Part \eqref{itemlabel:lem:Q:1}]
It follows from 
the property guaranteed by Part \eqref{itemlabel:coverproperty-of-H} of Lemma~\ref{lem:H}
that for each $\alpha \in I$, 
\[
C_\alpha \subset G^{+,o}(\bar\eps,\bar\delta)^c \cap[-R,R]^\ell,
\] 
\[
\bigcup_{\alpha \in I} C_\alpha = G^{+,o}(\bar\eps,\bar\delta)^c \cap[-R,R]^\ell
=
\bigcup_{\substack{\sigma \in \SIGN(\mathcal{H}), 
\\
\RR(\sigma) \subset \bigcup_{\alpha \in I} C_\alpha}}
\RR(\sigma,\rR\la\bar\eps,\bar\delta,\bar\zeta\ra^\ell).
\]

So it suffices to prove that $G^{+,o}(\bar\eps,\bar\delta) \cap [-R,R]^\ell \cap \rR^\ell = \emptyset$.
Since, $\lim_{\bar\delta} G^{+,o}(\bar\eps,\bar\delta) \cap [-R,R]^\ell = G(\bar\eps) \cap [-R,R]^\ell$,
in order to prove that $G^{+,o}(\bar\eps,\bar\delta) \cap [-R,R]^\ell \cap \rR^\ell = \emptyset$, it is enough to 
show that $G(\bar\eps) \cap [-R,R]^\ell \cap \rR^\ell = \emptyset$. 
Suppose there exists $\y \in G(\bar\eps) \cap [-R,R]^\ell \cap \rR^\ell$. Then there are two cases:
\begin{enumerate}[{Case }1.]
\item
There exists $\mathcal{Q} \subset \mathcal{P}^\star(\bar\eps)$, $0 < \card(\mathcal{Q}) \leq k$ such that
$\y \in \crit(\mathcal{Q})$.
But this would imply that there exists $(\x,\y) \in \rR\la\bar\eps\ra^{k+\ell}$ such that
the rank of the corresponding Jacobian matrix
\[
\left( \frac{\partial P}{\partial X_i}
\right)_{\substack{P \in \mathcal{Q} \\
1\leq i \leq k}}
\]
at the point $(\x,\y)$ is not full, which 
contradicts Lemma~\ref{lem:non-singular2}.

\item
There exists $\mathcal{Q} \subset \mathcal{P}^\star(\bar\eps)$, $k < \card(\mathcal{Q}) \leq k+\ell$ such that
$\y \in \pi_Y(\ZZ(\mathcal{Q},\rR\la\bar\eps\ra^{k+\ell}))$. But since $\card(\mathcal{Q}) > k$, and $\y \in \rR^\ell$,
by Lemma~\ref{lem:non-singular2}, 
$\ZZ(\mathcal{Q}(\cdot,\y), \rR\la\bar\eps\ra^k) = \emptyset$, which is a contradiction.
\end{enumerate}
\end{proof}
\begin{proof}[Proof of Part \eqref{itemlabel:lem:Q:2}]
Let $\eta = (\bar\eps,\bar\delta,\bar\zeta)$.
For each $H \in \mathcal{H}$, write 
\[
H =\sum_\alpha H_\alpha \eta^\alpha,
\]
where $\alpha \in \mathbb{N}^{|\bar\eps| + |\bar\delta| + |\bar\zeta|}$ is a multi-index exponent, and each $H_\alpha \in \rR[Y_1,\ldots,Y_\ell]$. 
Denote by  $\supp(H) = \{\alpha \mid| \mid H_\alpha \neq 0\}$.
It follows from  the bounds on the set $\mathcal{H}$ 
stated in 
Parts \eqref{itemlabel:complexity-of-H} and \eqref{itemlabel:dependence-of-H} of 
Lemma~\ref{lem:H}
that
$\card(\supp(H)) \leq (s d)^{(k +\ell)^{O(1)}}$.
Let 
\[
\mathcal{Q} = \bigcup_{H \in \mathcal{H}} \{H_\alpha \mid H \in \supp(H)\}.
\]

It follows from the definition of the ordering of the real closed field $\rR\la\eta\ra$ (cf. Notation~\ref{not:Puiseux}), that
for any $\y \in \R^\ell$, and $H \in \mathcal{H}$, $\sign(H(\y))$ is determined by the tuple of signs $(\sign(H_\alpha(\y)))_{\alpha \in \supp(H)}$.
It is now easy to see that $\mathcal{Q}$ satisfies the claimed properties. 
\end{proof}
\end{proof}

\begin{lemma}
\label{lem:homotopic}
Let $\y \in  \rR^\ell$. Then, there exists a semi-algebraic deformation retraction
$
S^\star(\bar\eps)_\y \rightarrow \E(S_\y, \rR\la\bar\eps\ra).
$
In particular,
$\E(S_\y, \rR\la\bar\eps\ra)$ 
is semi-algebraically homotopy equivalent to $S^\star(\bar\eps)_\y$.
\end{lemma}

\begin{proof}
This is a consequence of Lemma 16.17 in \cite{BPRbook2} noting that the set $S^\star(\bar\eps)_\y$ is bounded
over $\rR$.
\end{proof}

\begin{proof}[Proof of Proposition~\ref{prop:BElim}]
Let $\mathcal{Q}$ be as in Lemma~\ref{lem:Q}, and
let $\sigma \in \SIGN(\mathcal{Q})$,  such that $\RR(\sigma,\rR^\ell) \subset C_\alpha$ for some $\alpha \in I$
(following the same notation as in Lemma~\ref{lem:Q}).
Let $C \in \Cc(\RR(\sigma,\rR^\ell))$.

Now, let $D_1^\alpha, \ldots, D_N^\alpha \subset \rR\la\bar\eps,\bar\delta\ra^\ell$ 
be the  elements of $\Cc(D_\alpha)$ where 
$D_\alpha = \pi_Y^{-1}(C_\alpha) \cap \E(S^\star(\bar\eps) \cap [-R,R]^{k+\ell},\rR\la\bar\eps,\bar\delta\ra)$.

The proposition will follow from the following two claims.

\begin{claim}
\label{claim:1}
\[
\left(\lim_{\bar\eps,\bar\delta} D_\alpha\right)_C  =  S_C.
\]
\end{claim}

\begin{claim}
\label{claim:2}
\[
\left(\lim_{\bar\eps,\bar\delta}\left( D_\alpha^i\right)\right)_C, i =1,\ldots, N,
\]
are the semi-algebraically connected components of $S_C$.
\end{claim}

\begin{proof}[Proof of Claim~\ref{claim:1}]
First observe that since   $S^{\star,c}(\bar\eps) \cap [-R,R]^{k+\ell}$ is bounded over $\rR$, so is $D_\alpha$. Also, it is clear that 
\[
S_C \subset \left(\lim_{\bar\eps,\bar\delta} D_\alpha\right)_C.
\]
To prove the reverse inclusion, let 
$(\x,\y) \in \left(\lim_{\bar\eps,\bar\delta} D_\alpha\right)_C$. Then, there exists 
$(\tilde{\x},\tilde{\y}) \in   D_\alpha$, such that  
\[
\lim_{\bar\eps,\bar\delta} (\tilde{\x},\tilde{\y}) = (\x,\y).
\]
Since $S^{\star,c}(\bar\eps) \cap[-R,R]^{k+\ell}$  (cf. Eqn. \eqref{eqn:S-star}) is a closed semi-algebraic subset
of $\R\la\bar\eps\ra^{k+\ell}$ bounded over $\rR$, it follows that
\[
\lim_{\bar\delta} (E(S^{\star,c}(\bar\eps) \cap[-R,R]^{k+\ell}, \rR\la\bar\eps,\bar\delta\ra)) = 
S^{\star,c}(\bar\eps) \cap[-R,R]^{k+\ell}.
\]
This implies
\[
\lim_{\bar\delta} (\tilde{\x},\tilde{\y}) \in S^{\star,c}(\bar\eps)_{\lim_{\bar\delta} \tilde{\y}}.
\]

Finally, it follows from the fact
\[
\lim_{\bar\eps} \left(S^{\star,c} \cap [-R,R]^{k+\ell}\right) = S,
\]
that
\[
\lim_{\bar\eps,\bar\delta} (\tilde{\x},\tilde{\y}) \in S_{\y} \subset S_C.
\]
\end{proof}

\begin{proof}[Proof of Claim~\ref{claim:2}]
Using Claim~\ref{claim:1} it suffices to prove that for each $i, 1 \leq i \leq N$,
$\left(\lim_{\bar\eps,\bar\delta}\left( D_\alpha^i\right)\right)_C$ is semi-algebraically connected,
and for $1 \leq i < j \leq N$,
 \[
 \left(\lim_{\bar\eps,\bar\delta}\left( D_\alpha^i\right)\right)_C  \cap \left(\lim_{\bar\eps,\bar\delta}\left( D_\alpha^j\right)\right)_C  = \emptyset.
 \]
 
In order to prove that $\left(\lim_{\bar\eps,\bar\delta}\left( D_\alpha^i\right)\right)_C $ is semi-algebraically connected, let \[
(\x,\y), (\x',\y') \in \left(\lim_{\bar\eps,\bar\delta}\left( D_\alpha^i\right)\right)_C.
\] 
 
 Let $\gamma:[0,1] \rightarrow C$ be a semi-algebraic path, with $\gamma(0) = \y, \gamma(1) = \y'$ (which exists since $C$ is semi-algebraically connected).
 
 First observe, that since $\y,\y' \in C \subset \rR^\ell$, it follows from Lemma~\ref{lem:homotopic}, that
 $\x,\x' \in D^i_\alpha$.
 
 Then, using Lemma~\ref{lem:trivial-fibration}, there exists a lift of $\E(\gamma,\rR\la\bar\eps,\bar\delta\ra)$
 to a semi-algebraic path,
 \[
 \widehat{\gamma}[0,1] \rightarrow D^i_\alpha \cap\pi_Y^{-1}(\E(C,\rR\la\bar\eps,\bar\delta\ra)), 
 \]
 \begin{eqnarray*}
 \widehat{\gamma}(0) &=&(\x,\y), \\
 \widehat{\gamma}(1) &=& (\x',\y'),\\
 \pi_Y(\widehat{\gamma}(t)) &=& \E(\gamma,\rR\la\bar\eps,\bar\delta\ra)(t), 0 \leq t \leq 1.
 \end{eqnarray*}
 
 Then, 
 \[
 \lim_{\bar\eps,\bar\delta} \circ \widehat{\gamma}: [0,1] \rightarrow \left(\lim_{\bar\eps,\bar\delta}\left( D_\alpha^i\right)\right)_C
 \] 
 is a semi-algebraic path 
 connecting $(\x,\y)$ to $(\x',\y')$ with image in 
 $\left(\lim_{\bar\eps,\bar\delta}\left( D_\alpha^i\right)\right)_C$, proving that 
 $\left(\lim_{\bar\eps,\bar\delta}\left( D_\alpha^i\right)\right)_C$ is semi-algebraically connected.

We now prove that
for $1 \leq i < j \leq N$,
 \[
 \left(\lim_{\bar\eps,\bar\delta}\left( D_\alpha^i\right)\right)_C \cap \left(\lim_{\bar\eps,\bar\delta}\left( D_\alpha^j\right)\right)_C = \emptyset.
 \]
 
Suppose that 
\[
(\x,\y) \in \left(\lim_{\bar\eps,\bar\delta}\left( D_\alpha^i\right)\right)_C \cap \left(\lim_{\bar\eps,\bar\delta}\left( D_\alpha^j\right)\right)_C.
\] 
Using Lemma~\ref{lem:homotopic} and the fact that
 $\y \in C \subset \rR^\ell$, this would implies that 
 $(\x,\y) \in D^i_\alpha \cap D^j_\alpha$. But this is impossible, since $D^i_\alpha$ and $D^j_\alpha$ are distinct
 semi-algebraically connected components of $D_\alpha$.
 \end{proof} 

The proposition now follows from Claim~\ref{claim:2}.

\end{proof}


\bibliographystyle{plain}
\bibliography{master}

\def\cprime{$'$} \def\cprime{$'$}
\begin{thebibliography}{10}

\bibitem{Aguilar}
Marcelo Aguilar, Samuel Gitler, and Carlos Prieto.
\newblock {\em Algebraic topology from a homotopical viewpoint}.
\newblock Universitext. Springer-Verlag, New York, 2002.
\newblock Translated from the Spanish by Stephen Bruce Sontz.

\bibitem{B99}
S.~Basu.
\newblock On bounding the {B}etti numbers and computing the {E}uler
  characteristic of semi-algebraic sets.
\newblock {\em Discrete Comput. Geom.}, 22(1):1--18, 1999.

\bibitem{BPR10}
S.~Basu, R.~Pollack, and M.-F. Roy.
\newblock Betti number bounds, applications and algorithms.
\newblock In {\em Current Trends in Combinatorial and Computational Geometry:
  Papers from the Special Program at MSRI}, volume~52 of {\em MSRI
  Publications}, pages 87--97. Cambridge University Press, 2005.

\bibitem{BPR02}
S.~Basu, R.~Pollack, and M.-F. Roy.
\newblock On the {B}etti numbers of sign conditions.
\newblock {\em Proc. Amer. Math. Soc.}, 133(4):965--974 (electronic), 2005.

\bibitem{BPRbook2}
S.~Basu, R.~Pollack, and M.-F. Roy.
\newblock {\em Algorithms in real algebraic geometry}, volume~10 of {\em
  Algorithms and Computation in Mathematics}.
\newblock Springer-Verlag, Berlin, 2006 (second edition).

\bibitem{BV06}
S.~Basu and N.~Vorobjov.
\newblock On the number of homotopy types of fibres of a definable map.
\newblock {\em J. Lond. Math. Soc. (2)}, 76(3):757--776, 2007.

\bibitem{Basu-sheaf}
Saugata Basu.
\newblock A {C}omplexity {T}heory of {C}onstructible {F}unctions and {S}heaves.
\newblock {\em Found. Comput. Math.}, 15(1):199--279, 2015.

\bibitem{BCR}
J.~Bochnak, M.~Coste, and M.-F. Roy.
\newblock {\em G\'eom\'etrie alg\'ebrique r\'eelle (Second edition in english:
  Real Algebraic Geometry)}, volume 12 (36) of {\em Ergebnisse der Mathematik
  und ihrer Grenzgebiete [Results in Mathematics and Related Areas ]}.
\newblock Springer-Verlag, Berlin, 1987 (1998).

\bibitem{Rham}
O~Burlet and G~de~Rham.
\newblock Sur certaines applications g{\'e}n{\'e}riques d'une vari{\'e}t{\'e}
  sur certaines applications g{\'e}n{\'e}riques d'une vari{\'e}t{\'e} close
  {\`a} 2 dimensions dans le plan.
\newblock {\em L'Ensiegnement Math{\'e}matique}, 20:275--292, 1974.

\bibitem{Michel2}
Michel Coste.
\newblock {\em An introduction to o-minimal geometry}.
\newblock Istituti Editoriali e Poligrafici Internazionali, Pisa, 2000.
\newblock Dip. Mat. Univ. Pisa, Dottorato di Ricerca in Matematica.

\bibitem{de-Silva-et-al}
Vin de~Silva, Elizabeth Munch, and Amit Patel.
\newblock Categorified {R}eeb graphs.
\newblock {\em Discrete Comput. Geom.}, 55(4):854--906, 2016.

\bibitem{Dey-et-al-2017}
T.~K. {Dey}, F.~{Memoli}, and Y.~{Wang}.
\newblock {Topological Analysis of Nerves, Reeb Spaces, Mappers, and Multiscale
  Mappers}.
\newblock {\em ArXiv e-prints}, March 2017.

\bibitem{RSP}
Herbert Edelsbrunner, John Harer, and Amit~K. Patel.
\newblock Reeb spaces of piecewise linear mappings.
\newblock In {\em Proceedings of the Twenty-fourth Annual Symposium on
  Computational Geometry}, SCG '08, pages 242--250, New York, NY, USA, 2008.
  ACM.

\bibitem{Edelsbrunner-Harer}
Herbert Edelsbrunner and John~L. Harer.
\newblock {\em Computational topology}.
\newblock American Mathematical Society, Providence, RI, 2010.
\newblock An introduction.

\bibitem{GaV}
A.~Gabrielov and N.~Vorobjov.
\newblock Betti numbers of semialgebraic sets defined by quantifier-free
  formulae.
\newblock {\em Discrete Comput. Geom.}, 33(3):395--401, 2005.

\bibitem{GVZ04}
A.~Gabrielov, N.~Vorobjov, and T.~Zell.
\newblock Betti numbers of semialgebraic and sub-{P}faffian sets.
\newblock {\em J. London Math. Soc. (2)}, 69(1):27--43, 2004.

\bibitem{GV07}
Andrei Gabrielov and Nicolai Vorobjov.
\newblock Approximation of definable sets by compact families, and upper bounds
  on homotopy and homology.
\newblock {\em J. Lond. Math. Soc. (2)}, 80(1):35--54, 2009.

\bibitem{Godement}
Roger Godement.
\newblock {\em Topologie alg\'ebrique et th\'eorie des faisceaux}.
\newblock Actualit'es Sci. Ind. No. 1252. Publ. Math. Univ. Strasbourg. No. 13.
  Hermann, Paris, 1958.

\bibitem{Hartshorne}
Robin Hartshorne.
\newblock {\em Algebraic geometry}.
\newblock Springer-Verlag, New York-Heidelberg, 1977.
\newblock Graduate Texts in Mathematics, No. 52.

\bibitem{Milnor2}
J.~Milnor.
\newblock On the {B}etti numbers of real varieties.
\newblock {\em Proc. Amer. Math. Soc.}, 15:275--280, 1964.

\bibitem{Mimura-Toda-book}
Mamoru Mimura and Hirosi Toda.
\newblock {\em Topology of {L}ie groups. {I}, {II}}, volume~91 of {\em
  Translations of Mathematical Monographs}.
\newblock American Mathematical Society, Providence, RI, 1991.
\newblock Translated from the 1978 Japanese edition by the authors.

\bibitem{Munch}
E.~{Munch} and B.~{Wang}.
\newblock {Convergence between Categorical Representations of Reeb Space and
  Mapper}.
\newblock {\em ArXiv e-prints}, December 2015.

\bibitem{Patel}
Amit Patel.
\newblock {\em Reeb Spaces and the Robustness of Preimages}.
\newblock PhD thesis, Duke University, 2010.

\bibitem{OP}
I.~G. Petrovski{\u\i} and O.~A. Ole{\u\i}nik.
\newblock On the topology of real algebraic surfaces.
\newblock {\em Izvestiya Akad. Nauk SSSR. Ser. Mat.}, 13:389--402, 1949.

\bibitem{PS1}
A.~Pillay and C.~Steinhorn.
\newblock Definable sets in ordered structures. {I}.
\newblock {\em Trans. Amer. Math. Soc.}, 295(2):565--592, 1986.

\bibitem{PS2}
A.~Pillay and C.~Steinhorn.
\newblock Definable sets in ordered structures. {III}.
\newblock {\em Trans. Amer. Math. Soc.}, 309(2):469--576, 1988.

\bibitem{Reeb}
Georges Reeb.
\newblock Sur les points singuliers d'une forme de pfaff compl{\`e}tement
  int{\'e}grable ou d'une fonction num{\'e}rique.
\newblock {\em Comptes Rendus de l'Acad{\'e}mie des Sciences}, 222:847--849,
  1946.

\bibitem{Rolin}
J.-P. Rolin, P.~Speissegger, and A.~J. Wilkie.
\newblock Quasianalytic {D}enjoy-{C}arleman classes and o-minimality.
\newblock {\em J. Amer. Math. Soc.}, 16(4):751--777 (electronic), 2003.

\bibitem{Singh}
Gurjeet Singh, Facundo M{\'e}moli, and Gunnar Carlsson.
\newblock Topological methods for the analysis of high dimensional data sets
  and 3d object recognition.
\newblock {\em Eurographics Symposium of Point-Based Graphics}, 2007.

\bibitem{T}
R.~Thom.
\newblock Sur l'homologie des vari\'et\'es alg\'ebriques r\'eelles.
\newblock In {\em Differential and Combinatorial Topology (A Symposium in Honor
  of Marston Morse)}, pages 255--265. Princeton Univ. Press, Princeton, N.J.,
  1965.

\bibitem{Dries5}
L.~van~den Dries.
\newblock Remarks on {T}arski's problem concerning {$({\bf R},\,+,\,\cdot
  ,\,{\rm exp})$}.
\newblock In {\em Logic colloquium '82 (Florence, 1982)}, volume 112 of {\em
  Stud. Logic Found. Math.}, pages 97--121. North-Holland, Amsterdam, 1984.

\bibitem{Dries}
L.~van~den Dries.
\newblock {\em Tame topology and o-minimal structures}, volume 248 of {\em
  London Mathematical Society Lecture Note Series}.
\newblock Cambridge University Press, Cambridge, 1998.

\bibitem{Dries2}
L.~van~den Dries and C.~Miller.
\newblock Geometric categories and o-minimal structures.
\newblock {\em Duke Math. J.}, 84(2):497--540, 1996.

\bibitem{Dries3}
L.~van~den Dries and P.~Speissegger.
\newblock The real field with convergent generalized power series.
\newblock {\em Trans. Amer. Math. Soc.}, 350(11):4377--4421, 1998.

\bibitem{Dries4}
L.~van~den Dries and P.~Speissegger.
\newblock The field of reals with multisummable series and the exponential
  function.
\newblock {\em Proc. London Math. Soc. (3)}, 81(3):513--565, 2000.

\bibitem{Wilkie}
A.~J. Wilkie.
\newblock Model completeness results for expansions of the ordered field of
  real numbers by restricted {P}faffian functions and the exponential function.
\newblock {\em J. Amer. Math. Soc.}, 9(4):1051--1094, 1996.

\bibitem{Wilkie2}
A.~J. Wilkie.
\newblock A theorem of the complement and some new o-minimal structures.
\newblock {\em Selecta Math. (N.S.)}, 5(4):397--421, 1999.

\bibitem{Wilkie2009}
Alex~J. Wilkie.
\newblock o-minimal structures.
\newblock {\em Ast\'erisque}, (326):Exp. No. 985, vii, 131--142 (2010), 2009.
\newblock S\'eminaire Bourbaki. Vol. 2007/2008.

\end{thebibliography}

\end{document}